\documentclass[conference]{IEEEtran}

\usepackage[english]{babel}
\usepackage{bbm} 
\usepackage{url}

\usepackage{amsmath}
\usepackage{amsthm}
\usepackage{amssymb}

\usepackage{marvosym} 

\newcommand{\e}[1]{ {\mathrm{e}}^{ #1 } }

\newcommand{\iterand}[2]{ #1^{[#2]} }
\newcommand{\vect}[1]{ \boldsymbol{#1} }
\newcommand{\vectones}[1]{ \vect{1}_{#1} }

\newcommand{\vectComponent}[2]{ #1_{#2} }
\newcommand{\vectInLine}[1]{ ( #1 )^{\mathrm{T}} }
\newcommand{\matrixElement}[3]{ {#1}_{#2,#3} }

\newcommand{\pnorm}[2]{ \| #1 \|{}_{#2} }
\newcommand{\bigpnorm}[2]{ \Bigl\| #1 \Bigr\|_{#2} }
\newcommand{\transpose}[1]{ #1{}^{\mathrm{T}} }
\newcommand{\truncate}[2]{ [ {#1} ]^{#2} }
\newcommand{\bigtruncate}[2]{ \Bigl[ {#1} \Bigr]^{#2} }
\newcommand{\gradientOperatorWrt}[1]{ \nabla_{#1} }

\newcommand{\criticalpoint}[1]{  #1^{\textnormal{opt}} }

\newcommand{\cardinality}[1]{ | #1 | }

\newcommand{\expectation}[1]{ \mathbb{E} [ #1 ] }

\newcommand{\indicator}[1]{ \mathbbm{1} [ #1 ] }

\newcommand{\probability}[1]{ \mathbb{P} [ #1 ] }
\newcommand{\probabilityBig}[1]{ \mathbb{P} \Bigl[ #1 \Bigr] }

\newcommand{\process}[2]{ \{ #1 \}_{ #2 } }
\newcommand{\totalVariation}[1]{ || #1 ||_{\mathrm{var}} }

\newcommand{\varianceWrt}[2]{ \mathrm{Var}_{#2} [ #1 ] }

\newcommand{\naturalNumbers}{ \mathbb{N} }
\newcommand{\realNumbers}{ \mathbb{R} }
\newcommand{\positiveRealNumbers}{ [0,\infty) }
\newcommand{\strictlyPositiveRealNumbers}{ (0,\infty) }

\newcommand{\bigO}[1]{ \mathcal{O}(#1) }

\newcommand{\sA}{ {\ensuremath{x}} } 
\newcommand{\sB}{ {\ensuremath{y}} } 
\newcommand{\sC}{ {\ensuremath{z}} } 

\newcommand{\svA}{ {\sA} }
\newcommand{\svB}{ {\sB} }
\newcommand{\svC}{ {\sC} }

\newcommand{\stateSpace}{ \Omega }
\newcommand{\itstep}[1]{\iterand{a}{#1}}
\newcommand{\iterror}[1]{\iterand{e}{#1}}
\newcommand{\itfreq}[1]{\iterand{f}{#1}}
\newcommand{\totaltime}[1]{\iterand{t}{#1}}
\newcommand{\ErlangB}[2]{{B}(#1,#2)}
\newcommand{\ErlangBestimate}{\hat{{B}}}
\newcommand{\ExpectedB}[2]{{L}(#1,#2)}
\newcommand{\ExpectedBestimate}{\hat{{L}}}

\newcommand{\ErlangBwrt}[1]{{B}_{#1}}
\newcommand{\ExpectedBwrt}[1]{{L}_{#1}}

\newcommand{\infinitesimalGenerator}[1]{{ \mathcal{L} {#1}}}
\newcommand{\normalizationConstant}{Z}
\newcommand{\numberOfServers}{s}

\newcommand{\LipschitzConstant}{c_{\mathrm{l}}}
\newcommand{\ergodicConstant}{c_{\mathrm{e}}}
\newcommand{\gradientConstant}{c_{\mathrm{g}}}

\newcommand{\martingaleConstant}{c_{\mathrm{m}}}
\newcommand{\startingErrorConstant}{c_{\mathrm{r}}}

\newcommand{\QuodEratDemonstrandum}{\hfill \ensuremath{\Box}}

\newcommand{\refEquation}[1]{{\textrm{\eqref{#1}}}}

\newcommand{\refTheorem}[1]{{\textrm{Theorem~\ref{#1}}}}

\newcommand{\refProposition}[1]{{\textrm{Proposition~\ref{#1}}}}
\newcommand{\refLemma}[1]{{\textrm{Lemma~\ref{#1}}}}

\newcommand{\refSection}[1]{{\textrm{\S\ref{#1}}}}
\newcommand{\refAppendixSection}[1]{{\textrm{Appendix \ref{#1}}}}

\newtheorem{theorem}{Theorem}
\newtheorem{proposition}{Proposition}
\newtheorem{lemma}{Lemma}

\newtheorem{definition}{Definition}

\begin{document}

%
%

\title{Online Optimization of Product-Form Networks}

\author{\IEEEauthorblockN{Jaron Sanders}
\IEEEauthorblockA{Dept. of Math. \& Comp. Science\\
Eindhoven University of Technology\\
5612 AZ, Den Dolech 2\\
Eindhoven, The Netherlands\\
jaron.sanders@tue.nl}
\and
\IEEEauthorblockN{Sem C. Borst}
\IEEEauthorblockA{Dept. of Math. \& Comp. Science\\
Eindhoven University of Technology\\
5612 AZ, Den Dolech 2\\
Eindhoven, The Netherlands\\
s.c.borst@tue.nl}
\and
\IEEEauthorblockN{Johan S.H. van Leeuwaarden}
\IEEEauthorblockA{Dept. of Math. \& Comp. Science\\
Eindhoven University of Technology\\
5612 AZ, Den Dolech 2\\
Eindhoven, The Netherlands\\
j.s.h.v.leeuwaarden@tue.nl}}

\maketitle

\begin{abstract}
We develop an online gradient algorithm for optimizing the performance of product-form networks through online adjustment of control parameters. The use of standard algorithms for finding optimal parameter settings is hampered by the prohibitive computational burden of calculating the gradient in terms of the stationary probabilities. The proposed approach instead relies on measuring empirical frequencies of the various states through simulation or online operation so as to obtain estimates for the gradient. Besides the reduction in computational effort, a further benefit of the online operation lies in the natural adaptation to slow variations in ambient parameters as commonly occurring in dynamic environments. On the downside, the measurements result in inherently noisy and biased estimates. We exploit mixing time results in order to overcome the impact of the bias and establish sufficient conditions for convergence to a globally optimal solution.
\end{abstract}


\section{Introduction}

Markov processes provide a versatile framework for modelling a wide
variety of stochastic systems, ranging from communication networks
and data center applications to content dissemination systems and physical
or social interaction processes~\cite{Bremaud99,Kelly85,Liggett85}.
In particular, key performance measures of the system under consideration, e.g.~buffer occupancies, response times,
loss probabilities or user throughputs, can typically be expressed
in terms of the stationary distribution~$\vect{\pi}$ of the Markov process.

In many applications, the stationary distribution~$\vect{\pi}$, and hence
the performance measures or statistical properties, crucially depend
on system parameters~$\vect{r}$ that can be controlled, e.g.~admission
thresholds, service rates, link weights or resource capacities.
In those cases, the interest is often not so much in evaluating the
performance of the system for given parameter values, but rather
in finding parameter settings $\criticalpoint{\vect{r}}$ that optimize
the performance or achieve an optimal trade-off between service level
and costs.

Specifically, let $\bar{u}(\vect{\pi}(\vect{r}))$ be a function expressing the performance
objective (to be minimized) in terms of the stationary distribution
$\vect{\pi}(\vect{r})$ as function of the system parameters $\vect{r}$
and let $c(\vect{r})$ be a function representing possible cost associated
with $\vect{r}$, e.g.~capital expense or power consumption.
Introducing $u(\vect{r}) = \bar{u}(\vect{\pi}(\vect{r})) + c(\vect{r})$,
the problem of interest may then be mathematically formulated as finding
\begin{align}
\criticalpoint{\vect{r}} = \arg\min_{\vect{r}} u(\vect{r}).
\label{eqn:Introduction_minimization_problem}
\end{align}
It is worth observing here that the problem formulation differs from
the typical Markov decision processes \cite{Bertsekas96,Puterman94},
which focus on selecting optimal actions in various states rather than
identifying optimal parameter values.

Optimization problem \refEquation{eqn:Introduction_minimization_problem}
could in principle be solved using mathematical programming approaches
such as gradient-based schemes.
In addition to the usual convexity issues, however, a further
difficulty arises from the fact that the stationary distribution
$\vect{\pi}(\vect{r})$ is only implicitly determined as a function
of $\vect{r}$ by the balance equations and is rarely available in
explicit form, which severely complicates both the evaluation of the
objective function $u(\vect{r})$ and calculation of its gradient $\gradientOperatorWrt{\vect{r}} u(\vect{r})$.

In the present paper we develop a gradient approach to solve the
optimization problem~(\ref{eqn:Introduction_minimization_problem})
for a class of Markov processes with product-form distributions.
This class of processes arises in a rich family of stochastic models,
such as loss networks \cite{JLSSS08,Kelly91},
open and closed queueing networks \cite{BCMP79,Kelly79},
wireless random-access networks \cite{BKMS87,WK05}
and various types of interacting-particle systems
\cite{Bremaud99,Liggett85}.

As we will show, the partial derivatives 
$\partial \vect{\pi}(\vect{r}) / \partial \vect{r}$ for this class of
processes can be written as linear combinations of products of
stationary probabilities~$\vect{\pi}(\vect{r})$, thus reducing the computation of the gradient to the evaluation of the equilibrium distribution.
The problem that yet remains in many situations is that the stationary
probabilities involve a normalization constant whose calculation is
computationally intensive and potentially NP-hard~\cite{LMK94}.
This issue is particularly pertinent in the context of iterative
optimization algorithms such as gradient-based schemes, where partial
derivatives need to be calculated repeatedly.

In order to circumvent the computational burden of calculating
the stationary probabilities, we adopt a gradient approach which
relies on measuring the empirical frequencies of the various
states so as to estimate the partial derivatives.
Specifically, in each iteration we observe the stochastic process for
some time period through simulation or online operation,
and we then calculate estimates for the gradient based on the measured
time fractions of the various states. Although the number of states may be extremely large, it turns out that in many situations one only needs to track the time fractions
of aggregate states rather than all individual states, and that these
aggregate states can be observed in an entirely distributed fashion.
Besides the reduction in computational effort, a further benefit of
the online operation lies in the fact that the algorithm will
automatically adapt to slow variations in ambient parameters which
are fairly common in dynamic environments.

While the measurements bypass the computational effort of calculating
the stationary probabilities, they result in inherently {\em noisy\/}
and {\em biased\/} estimates for the gradient.
The issue of noisy estimates is paramount in the field of
stochastic approximation, where years of research have resulted in many
robust stochastic approximation schemes which can cope with various
stochastic processes and forms of random noise \cite{Borkar08,KY03}.
In contrast, biased estimates present a much trickier issue,
which is usually not accounted for in stochastic approximation schemes.
In order to neutralize the impact of the bias, we focus the attention
on the family of reversible processes within the above-mentioned
class of Markov processes with product-form distributions~\cite{Kelly79}.
For reversible processes, powerful results are known for mixing times
\cite{DS91,LPW08}, which allow us to derive sufficient conditions
guaranteeing convergence to the optimal solution
of~(\ref{eqn:Introduction_minimization_problem}).
Intuitively, the mixing times provide an indication for the period of
time that we need to observe the stochastic process in order to
overcome the impact of the bias.

As a further condition to ensure convergence to the globally optimal
solution of~(\ref{eqn:Introduction_minimization_problem}) rather than
a possible local optimum, we assume the optimization
objective $u(\vect{r})$ to be convex in $\vect{r}$.
While convexity is generally non-trivial to establish, this can be
easily verified for the broad class of so-called log-likelihood functions
\begin{align}
u(\vect{r}) = \bar{u}(\vect{\pi}(\vect{r})) 
= - \transpose{ \vect{\alpha} } \ln \vect{\pi}(\vect{r})
= - \sum_{\svA \in \stateSpace} \vectComponent{\alpha}{\svA}
\ln{\vectComponent{\pi}{\svA}}(\vect{r}),
\label{eqn:Introduction_loglikelihood_function}
\end{align}
where $\stateSpace$ denotes the state space of the process,
$\vectComponent{\alpha}{\svA}$ are fixed coefficients
and $\vectComponent{\pi}{\svA}(\vect{r})$ is the stationary probability
of state $\svA$.
Taking partial derivatives of
\refEquation{eqn:Introduction_loglikelihood_function}, we find that the
first-order conditions reduce to linear constraints in terms of the
stationary probabilities.
In other words, the problem of attaining target values for expectations of
functionals of the stationary distribution can be cast as an optimization
objective of the form \refEquation{eqn:Introduction_loglikelihood_function}. A special case of \refEquation{eqn:Introduction_loglikelihood_function} was recently investigated by Jiang and Walrand \cite{JiangWalrand2009,JiangWalrand2010}. Their goal was to achieve target throughput values in CSMA networks by using an algorithm that adjusts the access or backoff parameters (represented by the vector $\vect{r}$ in \refEquation{eqn:Introduction_loglikelihood_function}) using empirical arrival and service rates. This in fact provided valuable inspiration for the work presented here, where we extend the scope of such algorithms to general product-form Markov processes and a larger class of objective functions. These generalizations require a different approach to deal with the impact of bias, as discussed in \refSection{sec:Error_bias}. 

Further important related work is done by Marbach and Tsitsiklis \cite{MarbachTsitsiklis2001,MarbachTsitsiklis2003}, see also \cite{Cao2007} for further background. In \cite{MarbachTsitsiklis2001,MarbachTsitsiklis2003}, an algorithm similar in spirit to ours is considered - an algorithm that aims to tackle a parameter optimization problem by relying on measurement-based evaluation of a gradient. Their convergence proof also involves analysis of noisy and biased estimates and the generic use of Lyapunov functions and martingale arguments. However, their expression for the gradient is fundamentally different and hence the specific proof arguments substantially differ as well. Although \cite{MarbachTsitsiklis2001,MarbachTsitsiklis2003} can be applied to more general Markov processes and furnishes greater versatility in use, it does not take advantage of simplifications that arise from the specific structure of product-form distributions as in this paper. Most importantly, however, the algorithm in \cite{MarbachTsitsiklis2001,MarbachTsitsiklis2003} differs in its updating method, because it updates parameters whenever the process visits recurrent states. Knowing whether the entire system is in a recurrent state (and thus when to update) requires information about all components of the system, making the algorithm in \cite{MarbachTsitsiklis2001,MarbachTsitsiklis2003} global in nature. This differs from our algorithm and that presented in \cite{JiangWalrand2009,JiangWalrand2010}, which can be implemented in a distributed manner.

The remainder of the paper is organized as follows.
In \refSection{sec:Optimization_algorithm}, we present a detailed
problem formulation, develop our measurement-based optimization
algorithm and state our main results.
Some illustrative
application scenarios are described next in \refSection{sec:Achieving_expectation}.
In \refSection{sec:Proving_for_convergence}, we first identify
conditions in terms of the measurement noise and bias which ensure
the convergence of the algorithm, and we then prove that these
conditions are satisfied.

\section{Algorithm description}
\label{sec:Optimization_algorithm}

Throughout this paper, we denote by $\vectComponent{b}{i}$ the $i$-th component of vector $\vect{b}$. When taking a scalar function of an $n$-dimensional vector $\vect{b}$, we do this component-wise, i.e.~$\exp{ \vect{b} } = \vectInLine{ \exp \vectComponent{b}{1}, ..., \exp \vectComponent{b}{n} }$. If we have a $\cardinality{\stateSpace}$-dimensional vector $\vect{b}$ in which each component corresponds to some state $\svA \in \stateSpace$, we write $\vectComponent{b}{\svA}$ for that component of $\vect{b}$ that corresponds to state $\svA$.  Similarly, we denote by $\matrixElement{A}{i}{j}$ the element in row $i$, column $j$ of matrix $A$. If rows and/or columns correspond to states in $\stateSpace$, we write $\matrixElement{A}{\svA}{\svB}$ instead. Finally, we denote by $\vectones{n}$ the $n$-dimensional vector of which all components equal one.

\subsection{Gradient scheme}
\label{sec:Algorithm_decription__Gradient_scheme}

Consider a Markov process $\process{ X(t) }{t \geq 0}$ that is irreducible, reversible and has a finite state space $\stateSpace$. Let $\vect{\pi}(\vect{r})$ denote its steady-state probability vector as a function of $d$ parameters $\vect{r} = \vectInLine{ \vectComponent{r}{1}, ..., \vectComponent{r}{d} }$, which arises naturally if one has a closed-form expression for the stationary distribution. The most prominent examples are the product-form distributions
\begin{align}
\vect{\pi}(\vect{r}) 
= \frac{1}{\normalizationConstant(\vect{r})} \exp{ ( A \vect{r} + \vect{b} ) },
\label{eqn:Product_form__Equilibrium_distribution}
\end{align} 
where $A \in \realNumbers^{\cardinality{\stateSpace} \times d}$ is a matrix, $\vect{b} \in \realNumbers^{\cardinality{\stateSpace}}$ is a vector and $\normalizationConstant(\vect{r})$ is the normalization constant.

We consider the optimization problem
\begin{align}
\min_{\vect{r} \in \mathcal{R}} u( \vect{r} ), \label{eqn:Minimization_problem}
\end{align}
where $u( \vect{r} )$ denotes an objective function that we assume to be convex in $\vect{r}$ on a hypercube $\mathcal{R} \subset \realNumbers^d$, representing the feasible range for the parameters $\vect{r}$. We furthermore require that \refEquation{eqn:Minimization_problem} has a unique minimizer $\criticalpoint{\vect{r}} = \arg \min_{\vect{r} \in \mathcal{R}} u( \vect{r} )$, and we assume that the gradient of $u(\vect{r})$ can be written as a function of $\vect{\pi}(\vect{r})$ and $\vect{r}$, i.e.~$\gradientOperatorWrt{\vect{r}} u(\vect{r}) = \vect{g}(\vect{\pi}(\vect{r}), \vect{r})$ where $\gradientOperatorWrt{\vect{r}} = \vectInLine{ \partial / \partial \vectComponent{r}{1}, ..., \partial / \partial \vectComponent{r}{d} }$. For example when $c(\vect{r}) = 0$, the gradient of $u(\vect{r}) = \bar{u}(\vect{\pi}(\vect{r}))$ can be written as
\begin{align}
\frac{ \partial \bar{u}(\vect{\pi}(\vect{r})) }{ \partial \vectComponent{r}{i} }
= \sum_{\svA \in \stateSpace} \frac{ \partial \bar{u}(\vect{\pi}(\vect{r})) }{ \partial \vectComponent{\pi}{\svA}(\vect{r}) } \frac{ \partial \vectComponent{\pi}{\svA}(\vect{r}) }{ \partial \vectComponent{r}{i} }
\label{eqn:Product_form__Chain_rule}
\end{align}
for $i = 1, ..., d$. For the important case of product-form distributions in \refEquation{eqn:Product_form__Equilibrium_distribution}, we have
\begin{align}
\frac{ \partial \vectComponent{\pi}{\svA}(\vect{r}) }{ \partial \vectComponent{r}{i} }
&= \frac{1}{ \normalizationConstant(\vect{r})^2 } \Bigl( 
\normalizationConstant(\vect{r}) \matrixElement{A}{\svA}{i} \exp{ \vectComponent{ ( A \vect{r} + \vect{b} ) }{\svA} } \nonumber \\
&\phantom{=} - \exp{ \vectComponent{ ( A \vect{r} + \vect{b} ) }{\svA} } \sum_{\svB \in \stateSpace} \matrixElement{A}{\svB}{i} \exp{ \vectComponent{ ( A \vect{r} + \vect{b} ) }{\svB} } \Bigr) \nonumber \\
&= \vectComponent{\pi}{\svA}(\vect{r}) \Bigl( \matrixElement{A}{\svA}{i} - \sum_{ \svB \in \stateSpace } \matrixElement{A}{\svB}{i} \vectComponent{\pi}{\svB}(\vect{r}) \Bigr), \label{eqn:Product_form__Partial_derivative_to_r}
\end{align}
so that $\partial \bar{u}(\vect{\pi}(\vect{r})) / \partial \vectComponent{r}{i} = \vectComponent{g}{i}(\vect{\pi}(\vect{r}))$ and therefore $\gradientOperatorWrt{\vect{r}} u(\vect{r}) = \vect{g}(\vect{\pi}(\vect{r}))$. While for this example the gradient can be written as a function of only $\vect{\pi}(\vect{r})$, in \refSection{sec:Example__Optimizing_service_and_cost_trade_off} we will encounter an example for which it is more efficient to write the gradient as a function of both $\vect{\pi}(\vect{r})$ and $\vect{r}$. For a calculation of such partial derivatives in a more general case of product-form networks, we refer the reader to \cite{LiuNain1991}.

Our goal is to find $\criticalpoint{\vect{r}}$ and in order to do so, it is natural to consider the gradient algorithm
\begin{align}
\iterand{ \vect{r} }{n+1} 
= \truncate{ \iterand{ \vect{r} }{n} - \itstep{n+1} \iterand{ \vect{g} }{n+1} }{\mathcal{R}}, \label{eqn:Analytic_gradient_algorithm}
\end{align}
where $\iterand{ \vect{g} }{n+1} = \vect{g}( \vect{\pi}( \iterand{\vect{r}}{n} ), \iterand{\vect{r}}{n} )$, and $n \in \naturalNumbers$ indexes the iteration. The $\itstep{n} \in \strictlyPositiveRealNumbers$ denote the step sizes of the algorithm, and we define the truncation operator as follows.

\begin{definition}
\label{def:Hypercube_and_truncation}
For $\mathcal{R} \subset \realNumbers^d$ of the form
\begin{align}
\mathcal{R} = [ \vectComponent{\mathcal{R}}{1}^{\min}, \vectComponent{\mathcal{R}}{1}^{\max}] \times ... \times [\vectComponent{\mathcal{R}}{d}^{\min}, \vectComponent{\mathcal{R}}{d}^{\max}],
\end{align}
the truncation $\truncate{ \vect{r} }{ \mathcal{R} } \in \realNumbers^d$ of $\vect{r} \in \realNumbers^d$ is defined component-wise as
\begin{align}
\vectComponent{ \truncate{ \vect{r} }{ \mathcal{R} } }{i}
= \max \bigl\{ \vectComponent{\mathcal{R}}{i}^{\min}, \min \bigl\{ \vectComponent{\mathcal{R}}{i}^{\max}, \vectComponent{r}{i} \bigr\} \bigr\}.
\end{align}
\end{definition}

\subsection{Online gradient algorithm}
\label{sec:Algorithm_decription__Online_gradient_algorithm}

It is well known that under suitable assumptions on the objective function and step sizes, the gradient algorithm in \refEquation{eqn:Analytic_gradient_algorithm} generates a sequence $\iterand{\vect{r}}{n}$ that converges to the optimal solution $\criticalpoint{\vect{r}}$. We also come back to this at the end of \refSection{sec:Conditions_for_convergence}. Calculating the gradient, however, may be difficult in practice, because it depends on $\vect{\pi}(\vect{r})$, limiting the applicability of \refEquation{eqn:Analytic_gradient_algorithm}. 

Instead of using \refEquation{eqn:Analytic_gradient_algorithm}, we will estimate $\vect{\pi}(\vect{r})$ by observing the evolution of the system. These observations will take place during time intervals $[\totaltime{n}, \totaltime{n+1}]$, where $0 = \totaltime{0} < \totaltime{1} < ...$. At the end of each interval, say at time $\totaltime{n+1}$, our algorithm will change the current system parameters $\iterand{\vect{R}}{n}$ to new parameters $\iterand{\vect{R}}{n+1}$ based on its observations.

The stochastic process $\process{Y(t)}{t \geq 0}$ that describes the system is given by $Y(t) = \iterand{Z}{n}(t)$, where $n$ is such that $t \in [ \totaltime{n}, \totaltime{n+1} ]$. The process $\process{ \iterand{Z}{n}(t) }{ \totaltime{n} \leq t \leq \totaltime{n+1} }$ is a time-homogeneous Markov process, which starts in $\iterand{Z}{n-1}(\totaltime{n})$ and evolves according to the generator of $\process{ X(t) }{ t \geq 0 }$ that corresponds to parameters $\iterand{\vect{R}}{n}$. 

Let us now make precise how our algorithm observes the system and makes decisions. At time $\iterand{t}{n+1}$, marking the end of observation period $n+1$, we calculate
\begin{align}
\iterand{\vectComponent{\hat{\Pi}}{\svA}}{n+1}
&= \frac{1}{ \totaltime{n+1} - \totaltime{n} } \int_{ \totaltime{n} }^{ \totaltime{n+1} } \indicator{ \iterand{Z}{n}(t) = \svA } dt
\end{align}
for every state $\svA \in \stateSpace$. During each interval, one thus keeps track of the fractions of time that the system is in every state. This constitutes an empirical estimate of $\vect{\pi}(\iterand{\vect{R}}{n})$. We then estimate the gradient $\iterand{ \vect{G} }{n+1} = \vect{g}( \vect{\pi}(\iterand{\vect{R}}{n}), \iterand{\vect{R}}{n} )$ by $\iterand{ \vect{\hat{G}} }{n+1} = \vect{g}( \iterand{\vect{\hat{\Pi}}}{n+1}, \iterand{\vect{R}}{n} )$. If we then apply \refEquation{eqn:Analytic_gradient_algorithm} using the estimated gradient instead of the actual gradient, we are essentially using the stochastic gradient algorithm
\begin{align}
\iterand{ \vect{R} }{n+1} 
= \truncate{ \iterand{ \vect{R} }{n} - \itstep{n+1} \iterand{ \vect{\hat{G}} }{n+1} }{\mathcal{R}}
\label{eqn:Theorem_Convergence__Stochastic_gradient_algorithm}
\end{align}
to update the parameters. 

Note that algorithm \refEquation{eqn:Analytic_gradient_algorithm} is deterministic, whereas \refEquation{eqn:Theorem_Convergence__Stochastic_gradient_algorithm} is stochastic. Also note that because we are estimating the gradient instead of explicitly calculating it, the algorithm in \refEquation{eqn:Theorem_Convergence__Stochastic_gradient_algorithm} is no longer guaranteed to converge to $\criticalpoint{\vect{r}}$. 

\subsection{Main result}

We now present technical assumptions which will guarantee convergence of \refEquation{eqn:Theorem_Convergence__Stochastic_gradient_algorithm}. For this, we need an additional sequence $\iterror{n}$ which we shall refer to as the error. It is related to the maximum allowable error when estimating the steady-state probability vector, which will be made precise in \refSection{sec:Error_bias}.

We require the sequences $\itstep{n}$, $\iterror{n}$ and $\itfreq{n} = 1 / ( \totaltime{n} - \totaltime{n-1} )$ to be such that 
\begin{align}
&\sum_{n=1}^{\infty} \itstep{n} = \infty, 
\,\,\, \sum_{n=1}^{\infty} (\itstep{n})^2 < \infty, \label{eqn:Theorem_Convergence__Assumption_decreasing_stepsizes}
\end{align}
and
\begin{align}
&\sum_{n=1}^{\infty} \itstep{n} \iterror{n} < \infty, 
\,\,\, \sum_{n=1}^{\infty} \itstep{n} \exp{ \Bigl( - \frac{ (\iterror{n})^2 }{ 4 \cardinality{\stateSpace}^2 \kappa \itfreq{n} } \Bigr) } < \infty, \label{eqn:Theorem_Convergence__Assumption_decreasing_error_and_update_frequency}
\end{align}
for any $\kappa \in \strictlyPositiveRealNumbers$. We also require boundedness and regularity of $\vect{g}(\vect{\pi}(\vect{r}),\vect{r})$, in the sense that there exist constants $\gradientConstant, \LipschitzConstant \in \positiveRealNumbers$ such that
\begin{align}
| \vectComponent{g}{i}( \vect{\mu}, \vect{r} ) - \vectComponent{g}{i}( \vect{\nu}, \vect{r} ) |
&\leq \LipschitzConstant \totalVariation{ \vect{\mu} - \vect{\nu} } \textrm{ for } i = 1,...,d, \label{eqn:Theorem_Convergence__Assumption_Lipschitz_continuity} \\
\pnorm{ \vect{g}(\vect{\mu},\vect{r}) }{2} &\leq \gradientConstant,
\label{eqn:Theorem_Convergence__Assumption_Boundedness_of_g_in_r}
\end{align}
for all probability vectors $\vect{\mu}, \vect{\nu}$ and all $\vect{r} \in \mathcal{R}$. Here, $\totalVariation{ \vect{\mu} - \vect{\nu} } = \frac{1}{2} \sum_{\svA \in \stateSpace} | \vectComponent{\mu}{\svA} - \vectComponent{\nu}{\svA} |$ is the total variation distance. Under conditions \refEquation{eqn:Theorem_Convergence__Assumption_decreasing_stepsizes} -- \refEquation{eqn:Theorem_Convergence__Assumption_Boundedness_of_g_in_r} and the assumptions in \refSection{sec:Algorithm_decription__Gradient_scheme} and \refSection{sec:Algorithm_decription__Online_gradient_algorithm}, the following result holds.

\begin{theorem}
\label{thm:Convergence}
The sequence $\iterand{\vect{R}}{n}$ generated by the online algorithm \refEquation{eqn:Theorem_Convergence__Stochastic_gradient_algorithm} converges to the optimal solution $\criticalpoint{\vect{r}}$ of the optimization problem \refEquation{eqn:Minimization_problem} with probability one.
\end{theorem}

Condition \refEquation{eqn:Theorem_Convergence__Assumption_decreasing_stepsizes} is typical in stochastic approximation. It ensures that step sizes become smaller as $n$ increases, while remaining large enough so that the algorithm does not get stuck in a suboptimal solution. Condition \refEquation{eqn:Theorem_Convergence__Assumption_decreasing_error_and_update_frequency} then requires that the error $\iterror{n}$ for which we allow when estimating the steady-state probability vector must decrease. In order to guarantee this, the observation frequency $\itfreq{n}$ must eventually become smaller than the error, i.e.~$(\iterror{n})^2 / \itfreq{n} \rightarrow \infty$ as $n \rightarrow \infty$. Condition \refEquation{eqn:Theorem_Convergence__Assumption_Lipschitz_continuity} ensures that when we approximate the gradient of $u(\vect{r})$ by using empirical distributions that come increasingly closer to the actual $\vect{\pi}(\vect{r})$, our approximation of the gradient also comes increasingly closer to the actual gradient. It is the most non-trivial of all conditions and verification can be cumbersome. In \refSection{sec:Achieving_expectation} we discuss two illustrative examples for which \refEquation{eqn:Theorem_Convergence__Assumption_Lipschitz_continuity} holds. Lastly, condition \refEquation{eqn:Theorem_Convergence__Assumption_Boundedness_of_g_in_r} guarantees that the gradient does not explode, preventing the algorithm from making extremely large errors.

It is not difficult to define sequences that satisfy \refEquation{eqn:Theorem_Convergence__Assumption_decreasing_stepsizes} and \refEquation{eqn:Theorem_Convergence__Assumption_decreasing_error_and_update_frequency}. For example, setting $\itstep{n} = n^{-1}$, $\itfreq{n} = n^{-2 \alpha - \beta}$ and $\iterror{n} = n^{-\alpha}$ with $\alpha, \beta > 0$ suffices. In particular, note that for $\alpha = \beta = 1/3$, we have $\itstep{n} = n^{-1}$ and $\totaltime{n+1} - \totaltime{n} = n + 1$, which expresses that the algorithm should take smaller steps as time increases, while simultaneously lengthening the observation period.

The choices for $\itstep{n}$, $\iterror{n}$ and $\itfreq{n}$ strongly influence the behavior of the algorithm. Consider for instance the following two cases. Setting $\itstep{n} = n^{-1/2-\alpha}$ with $0 < \alpha \ll 1/2$ so that it barely satisfies \refEquation{eqn:Theorem_Convergence__Assumption_decreasing_stepsizes}, allows us to let $\iterror{n}$ decrease as slowly as $\iterror{n} = n^{-1/2}$. By \refEquation{eqn:Theorem_Convergence__Assumption_decreasing_error_and_update_frequency} we then need that $\itfreq{n} < n^{-1}$ or $\totaltime{n} - \totaltime{n-1} > n$. If we now consider the faster decreasing step size $\itstep{n} = n^{-1}$, which also barely satisfies \refEquation{eqn:Theorem_Convergence__Assumption_decreasing_stepsizes}, we find that a much slower decreasing $\iterror{n} = n^{-\alpha}$ with $0 < \alpha \ll 1$ suffices, implying by \refEquation{eqn:Theorem_Convergence__Assumption_decreasing_error_and_update_frequency} that $\itfreq{n} < n^{-2\alpha}$ or $\totaltime{n} - \totaltime{n-1} > n^{2\alpha}$ is required. From these two cases, one sees that smaller step sizes allow for shorter observation periods (recall that $0 < \alpha \ll 1$). The search for optimal settings of $\itstep{n}$, $\iterror{n}$ and $\itfreq{n}$ is an important topic for future research.

\section{Example applications}
\label{sec:Achieving_expectation}

We now discuss two example scenarios in which \refTheorem{thm:Convergence} can be applied. The first scenario concerns the optimal trade-off between performance and costs in an Erlang loss system. The second scenario considers a log-likelihood function as an objective function in combination with product-form stationary distributions. We should stress that these two examples, particularly the first one, primarily serve to illuminate the core features of our algorithm in relatively simple settings. These scenarios are not meant to reflect the full scope or unique realm of our algorithm and could conceivably also be tackled via alternative methods.

\subsection{Optimizing service, cost trade-off}
\label{sec:Example__Optimizing_service_and_cost_trade_off}

Consider the $M/M/\numberOfServers/\numberOfServers$ queue. Customers arrive according to a Poisson process with rate $\lambda$ and each customer has an exponentially distributed service requirement with unit mean. Each of the $\numberOfServers$ parallel servers works at rate $r$. The steady-state probability of $\sA \in \stateSpace = \{ 0, 1, ..., \numberOfServers \}$ customers in the system is then given by
\begin{align}
\vectComponent{\pi}{\sA}(r) = \frac{ \bigl( \lambda / r \bigr)^\sA / \sA! }{ \sum_{\sB=0}^{\numberOfServers} \bigl( \lambda / r \bigr)^\sB / \sB! }.
\end{align}
The steady-state probability that an arriving customer finds all servers occupied and is blocked is given by the Erlang loss formula $\ErlangB{\numberOfServers}{r} = \vectComponent{\pi}{\numberOfServers}(r)$. The mean stationary queue length is given by $\ExpectedB{\numberOfServers}{r} = \sum_{\sA = 1}^\numberOfServers \sA \vectComponent{\pi}{\sA}(r)$, and by Little's law, $\ExpectedB{\numberOfServers}{r} = \lambda ( 1 - \ErlangB{\numberOfServers}{r} ) / r$.

Suppose now that we want to minimize $\ErlangB{\numberOfServers}{r}$ by adjusting $r$ and that the costs of operating at service rate $r$ equal $c(r)$. Assume $c(r)$ to be convex in $r$ and its derivative $c'(r)$ to be bounded for all $r \in \mathcal{R}$. We thus aim to minimize $u(r) = \ErlangB{\numberOfServers}{r} + c(r)$. This objective function is convex in $r$ \cite{Harel1990}. Furthermore,
\begin{align}
g( \vect{\pi}(r), r ) 
&= \frac{ \ErlangB{\numberOfServers}{r} ( \ExpectedB{\numberOfServers}{r} - \numberOfServers ) }{ r } + c'(r), \label{eqn:Example_Erlang_B__Gradient}
\end{align}
for which we prove the following result in \refAppendixSection{sec:Appendix__Lipschitz_continuity}.

\begin{lemma}
\label{lemma:Lipschitz_continuity_of_g_for_MMss_queue}
If $\mathcal{R} = [\mathcal{R}^{\min}, \mathcal{R}^{\max}]$ with $0 < \mathcal{R}^{\min} < \mathcal{R}^{\max} < \infty$ and $g( \vect{\mu}, r )$ is given by \refEquation{eqn:Example_Erlang_B__Gradient}, then there exists constants $\gradientConstant, \LipschitzConstant \in \positiveRealNumbers$ such that
conditions \refEquation{eqn:Theorem_Convergence__Assumption_Lipschitz_continuity}, \refEquation{eqn:Theorem_Convergence__Assumption_Boundedness_of_g_in_r} hold for all probability vectors $\vect{\mu}, \vect{\nu}$ and all $r \in \mathcal{R}$.
\end{lemma}

Using \refLemma{lemma:Lipschitz_continuity_of_g_for_MMss_queue} we conclude that all conditions of \refTheorem{thm:Convergence} are met and that the gradient algorithm 
\begin{align}
\iterand{R}{n+1} = \bigtruncate{ \iterand{R}{n} - \itstep{n+1} \Bigl( \frac{ \iterand{ \ErlangBestimate }{n+1} ( \iterand{ \ExpectedBestimate }{n+1} - \numberOfServers ) }{ \iterand{R}{n} } + c'( \iterand{R}{n} ) \Bigr) }{\mathcal{R}} \nonumber
\end{align}
converges to the optimal solution. Here, $\iterand{ \ErlangBestimate }{n+1} = \iterand{ \vectComponent{\hat{\Pi}}{\numberOfServers} }{n+1}$ denotes an estimate of the loss probability and $\iterand{ \ExpectedBestimate }{n+1} = \sum_{\sA = 1}^{\numberOfServers} \sA \iterand{ \vectComponent{\hat{\Pi}}{\sA} }{n+1}$ denotes an estimate of the mean queue length.

\subsection{Log-likelihood and product forms}

Consider the log-likelihood function as defined in \refEquation{eqn:Introduction_loglikelihood_function} as objective function. We prove the following result in \refAppendixSection{sec:Appendix__Convexity_of_the_log_likelihood_function}.

\begin{lemma}
\label{lemma:Convexity_of_the_log_likelihood_function}
If $\vect{\pi}(\vect{r})$ satisfies the product form \refEquation{eqn:Product_form__Equilibrium_distribution}, then the log-likelihood function $u(\vect{r})$ in \refEquation{eqn:Introduction_loglikelihood_function} is convex in $\vect{r}$. 
\end{lemma}

Using $\partial \bar{u}(\vect{\pi}(\vect{r})) / \partial \vectComponent{\pi}{\svA} = - \vectComponent{\alpha}{\svA} / \vectComponent{\pi}{\svA}$ and substituting \refEquation{eqn:Product_form__Partial_derivative_to_r} into \refEquation{eqn:Product_form__Chain_rule} yields
\begin{align}
\vectComponent{g}{i}(\vect{\pi}(\vect{r}))
= \sum_{\svA \in \stateSpace} \vectComponent{\alpha}{\svA} \Bigl( \sum_{ \svB \in \stateSpace } \matrixElement{A}{\svB}{i} \vectComponent{\pi}{\svB}(\vect{r}) - \matrixElement{A}{\svA}{i} \Bigr). \label{eqn:Product_form__Gradient_of_log_likelihood_function_arbitrary_alpha}
\end{align}
We will only consider $\vect{\alpha} \in (0,1)^{\cardinality{\stateSpace}}$ that are probability vectors, so that $\transpose{ \vectones{\cardinality{\stateSpace}} } \vect{\alpha} = 1$. We can then interpret \refEquation{eqn:Product_form__Gradient_of_log_likelihood_function_arbitrary_alpha} as the difference between the expectation with respect to $\vect{\pi}(\vect{r})$, denoted by $\vectComponent{ ( \transpose{ A } \vect{\pi}(\vect{r}) ) }{i} = \sum_{ \svB \in \stateSpace } \matrixElement{A}{\svB}{i} \vectComponent{\pi}{\svB}(\vect{r})$, and the expectation with respect to $\vect{\alpha}$, denoted by $\vectComponent{ ( \transpose{ A } \vect{\alpha} ) }{i} = \sum_{ \svA \in \stateSpace } \matrixElement{A}{\svA}{i} \vectComponent{\alpha}{\svA}$, so that
\begin{align}
\vect{g}(\vect{\pi}(\vect{r}))
= \transpose{ A } \vect{\pi}(\vect{r}) - \transpose{ A } \vect{\alpha}.
\label{eqn:Product_form__Gradient_of_log_likelihood_function}
\end{align}
We assume that $\criticalpoint{\vect{r}}$ lies in the interior of $\mathcal{R}$, in which case optimality requires $\vect{g}(\vect{\pi}(\criticalpoint{\vect{r}})) = \vect{0}$ and thus $\transpose{ A } \vect{\pi}(\criticalpoint{\vect{r}}) = \transpose{ A } \vect{\alpha}$. We call $\vect{\gamma} = \transpose{ A } \vect{\alpha}$ the target vector, a name inspired by the fact that our algorithm seeks $\criticalpoint{\vect{r}}$ such that $\transpose{ A } \vect{\pi}(\criticalpoint{\vect{r}}) = \vect{\gamma}$.

Because $u(\vect{r})$ is convex in $\vect{r}$ and the target $\vect{\gamma}$ is achieved by the solution $\criticalpoint{\vect{r}}$ of \refEquation{eqn:Minimization_problem}, we want to use our online gradient algorithm \refEquation{eqn:Theorem_Convergence__Stochastic_gradient_algorithm} to find $\criticalpoint{\vect{r}}$. From \refEquation{eqn:Product_form__Gradient_of_log_likelihood_function}, it follows that
$| \vectComponent{g}{i}( \vect{\mu} ) - \vectComponent{g}{i}( \vect{\nu} )  \leq 2 \max_{\svA, i} \{ | \matrixElement{A}{\svA}{i} | \} \totalVariation{ \vect{\mu} - \vect{\nu} }$ for $i=1,...,d$, and that $\pnorm{ \vect{g}( \vect{\mu}, \vect{r} ) }{2} \leq \cardinality{\stateSpace} d \max_{x,i} \{ | \matrixElement{A}{x}{i} | \}$,
so that  \refEquation{eqn:Theorem_Convergence__Assumption_Lipschitz_continuity} and \refEquation{eqn:Theorem_Convergence__Assumption_Boundedness_of_g_in_r} are satisfied. Using \refTheorem{thm:Convergence}, we then arrive at the following result.

\begin{theorem}
\label{thm:Convergence__Product_form_steady_state_probability_vector}
Given any $\vect{\gamma} \in \realNumbers^{d}$ for which there exists an $\criticalpoint{\vect{r}}$ in the interior of $\mathcal{R}$ so that $\transpose{ A } \vect{\pi}(\criticalpoint{\vect{r}}) = \vect{\gamma}$, the online gradient algorithm 
\begin{align}
\iterand{ \vect{R} }{n+1} 
= \truncate{ \iterand{ \vect{R} }{n} - \itstep{n+1} \bigl( \transpose{ A } \iterand{ \vect{\hat{\Pi}} }{n+1} - \vect{\gamma} \bigr) }{\mathcal{R}} \label{eqn:Theorem_Convergence__Product_form_steady_state_probability_vectors}
\end{align}
converges to $\criticalpoint{\vect{r}}$ with probability one.
\end{theorem}

As an illustrative example, consider a loss network consisting of $L$~links
with capacities $\vect{c} = \vectInLine{ \vectComponent{c}{1}, ..., \vectComponent{c}{L} }$ shared by $K$ customer classes. Class-$k$ customers arrive according to a Poisson process with rate $\vectComponent{\lambda}{k}$ and require exponentially distributed holding times with mean $1 / \vectComponent{\mu}{k}$. Each class-$k$ customer requires capacity $\matrixElement{B}{k}{l}$ on link $l$ for the duration of its holding time, i.e.~$\matrixElement{B}{k}{l} = \vectComponent{b}{k} \matrixElement{J}{k}{l}$, where $\vectComponent{b}{k}$ is the nominal capacity requirement of a class-$k$ customer and $\matrixElement{J}{k}{l}$ has the value $0$ or $1$, indicating whether the route of class-$k$ customers contains link $l$ or not. When an arriving class-$k$ customer finds insufficient capacity available, it is blocked and lost. Denote the number of class-$k$ customers in the network at time~$t$ by $\vectComponent{X}{k}(t)$ and define $\vect{X}(t) = \vectInLine{ \vectComponent{X}{1}(t), ..., \vectComponent{X}{K}(t) }$. Under these assumptions, $\process{ \vect{X}(t) }{t \geq 0}$ is a reversible Markov process with state space $\stateSpace = \{ \vect{\svA} \in \naturalNumbers^K | B \vect{\svA} \leq \vect{c} \}$ and steady-state probability vector
\begin{align}
\vectComponent{\pi}{\vect{\svA}}(\vect{\rho}) = \frac{1}{\normalizationConstant(\vect{\rho})} \prod\limits_{k = 1}^{K} \frac{ ( \vectComponent{\rho}{k} )^{\vectComponent{\sA}{k}} }{ \vectComponent{\sA}{k}! }, 
\textrm{ where } 
\normalizationConstant(\vect{\rho}) = \sum_{\svB \in \stateSpace} \prod\limits_{k = 1}^{K} \frac{(\vectComponent{\rho}{k})^{ \vectComponent{\sB}{k} } }{ \vectComponent{\sB}{k}! }. \nonumber
\end{align}
Here, $\vectComponent{\rho}{k} = \vectComponent{\lambda}{k} / \vectComponent{\mu}{k}$ denotes the offered traffic of class $k$. Rewriting gives
\begin{align}
\vectComponent{\pi}{\vect{\svA}}(\vect{\rho})
= \frac{1}{\normalizationConstant(\vect{\rho})} \exp{ \Bigl( \sum_{k = 1}^{K} \vectComponent{\sA}{k} \ln \vectComponent{\rho}{k} - \ln ( \vectComponent{\sA}{k}! ) \Bigr) },
\end{align}
which matches \refEquation{eqn:Product_form__Equilibrium_distribution} with $d = K$, $\vectComponent{r}{k} = \ln \vectComponent{\rho}{k}$, $\matrixElement{A}{\svA}{k} = \vectComponent{\sA}{k}$ and $\vectComponent{b}{\svA} = - \sum_{k=1}^{K} \ln ( \vectComponent{\sA}{k}! )$. Note that $\vectComponent{ ( \transpose{ A } \vect{\pi}(\vect{r}) ) }{k} = \sum_{\svB \in \stateSpace} \vectComponent{\sB}{k} \vectComponent{\pi}{\sB}$ is the carried traffic of class $k$, i.e.~the steady-state average number of class-$k$ customers in the system, which we can empirically estimate by observing the system. We apply our algorithm by setting 
\begin{align}
\iterand{\vect{\rho}}{n+1} = \exp \Bigl( \truncate{ \ln \iterand{\vect{\rho}}{n} - \itstep{n+1} \bigl( \transpose{ A } \iterand{\vect{\hat{\Pi}}}{n+1} - \vect{\gamma} \bigr) }{\mathcal{R}} \Bigr), \label{eqn:Algorithm_for_log_likelihood_Poisson_class_system_example}
\end{align}
in order to adjust the amount of offered traffic $\vect{\rho}$ so as to achieve target carried traffic levels $\vect{\gamma}$. In practice, network operators usually have limited control over the amount of offered traffic, but they can typically adjust route selections fairly easily so as to achieve target blocking levels for a given offered traffic volume. Variations of the above algorithm can be used in such scenarios but go beyond the scope of the present paper.

In related work, Jiang and Walrand \cite{JiangWalrand2009,JiangWalrand2010} present an algorithm for achieving target throughputs in wireless CSMA networks. Their model can be interpreted as a special case of a loss network with unit link capacities. Their algorithm and convergence proof are therefore special cases of \refTheorem{thm:Convergence__Product_form_steady_state_probability_vector}.

\section{Convergence proof}
\label{sec:Proving_for_convergence}

We will now prove \refTheorem{thm:Convergence}. In \refSection{sec:Conditions_for_convergence}, we first explain our notion of convergence and then derive conditions on the error bias and zero-mean noise so that convergence is guaranteed. In \refSection{sec:Proving_convergence}, we show that under the assumptions of \refTheorem{thm:Convergence}, the error bias and zero-mean noise indeed satisfy the conditions derived in \refSection{sec:Conditions_for_convergence}. 

\subsection{Conditions for convergence}
\label{sec:Conditions_for_convergence}

\refTheorem{thm:Convergence} states that $\iterand{\vect{R}}{n}$ converges to $\criticalpoint{\vect{r}}$ with probability one. In order to prove that, we will establish that the following two properties hold for arbitrary $\delta, \varepsilon > 0$. As our first property, we want that $\iterand{\vect{R}}{n}$ comes close to $\criticalpoint{\vect{r}}$ infinitely often. We make this precise by requiring that for any $\delta > 0$, the set $\mathcal{H}_\delta = \{ \vect{r} \in \realNumbers^d | u(\vect{r}) \leq u( \criticalpoint{\vect{r}} ) + \delta / 2 \}$ is recurrent for $\process{ \iterand{\vect{R}}{n} }{ n \in \naturalNumbers }$. As our second property, we want that once $\iterand{\vect{R}}{n}$ comes close to $\criticalpoint{\vect{r}}$, it stays close to $\criticalpoint{\vect{r}}$ for all future iterations. Mathematically, we require that there exists an $m \in \naturalNumbers$ large enough so that $\pnorm{ \iterand{\vect{R}}{n} - \criticalpoint{\vect{r}} }{2}^2 \leq \pnorm{ \iterand{\vect{R}}{m} - \criticalpoint{\vect{r}} }{2}^2 + \varepsilon$ for all $n \geq m$, which we will call capture of $\iterand{\vect{R}}{n}$. 

We shall relate both recurrence and capture to the error bias and zero-mean noise, defined as $\iterand{\vect{B}}{n} = \expectation{ \iterand{ \vect{\hat{G}} }{n} | \iterand{\mathcal{F}}{n-1} } - \iterand{ \vect{G} }{n}$ and $\iterand{\vect{E}}{n} = \iterand{ \vect{\hat{G}} }{n} - \expectation{ \iterand{ \vect{\hat{G}} }{n} | \iterand{\mathcal{F}}{n-1} }$, respectively. Here, $\iterand{\mathcal{F}}{n-1}$ denotes the $\sigma$-field generated by the random vectors $\iterand{\vect{Z}}{0}, \iterand{\vect{Z}}{1}, ..., \iterand{\vect{Z}}{n-1}$, where $\iterand{\vect{Z}}{0} = \vectInLine{ \iterand{\vect{R}}{0}, X(0) }$ and $\iterand{\vect{Z}}{n} = \vectInLine{ \iterand{\vect{\hat{G}}}{n}, \iterand{\vect{R}}{n}, X({\totaltime{n}}) }$ for $n \geq 1$.

\subsubsection{Recurrence}

We begin with deriving conditions under which the set $\mathcal{H}_\delta$ is recurrent for $\process{ \iterand{\vect{R}}{n} }{ n \in \naturalNumbers }$, using the following result.

\begin{lemma}[\cite{KY03}, p.~115]
\label{lemma:KushnerYin2003__Recurrence_theorem}
Let $\process{ \iterand{\vect{R}}{n} }{ n }$ be an $\realNumbers^d$-valued stochastic process, not necessarily a Markov process. Let $\{ \iterand{\mathcal{F}}{n} \}$ be a sequence of nondecreasing $\sigma$-algebras, with $\iterand{\mathcal{F}}{n}$ measuring at least $\{ \iterand{\vect{R}}{i} | i \leq n \}$. Assume that $\itstep{n+1}$ are positive $\iterand{\mathcal{F}}{n}$-measurable random variables tending to zero with probability one and $\sum_{n} \itstep{n} = \infty$ with probability one. Let $V(\vect{r}) \geq 0$ and suppose that there are $\delta > 0$ and compact $\mathcal{H}_\delta \subset \realNumbers^d$ such that for all large $n$ and all $\vect{r} \not\in \mathcal{H}_\delta$,
\begin{align}
\expectation{ V(\iterand{\vect{R}}{n+1}) | \iterand{\mathcal{F}}{n} } - V(\iterand{\vect{R}}{n}) \leq - \itstep{n+1} \delta < 0. \label{eqn:KushnerYin2003__Inequality_to_guarantee_recurrence}
\end{align}
Then the set $\mathcal{H}_\delta$ is recurrent for $\process{\iterand{\vect{R}}{n}}{n \geq 0}$ in the sense that $\iterand{\vect{R}}{n} \in \mathcal{H}_\delta$ for infinitely many $n$ with probability one.
\end{lemma}

Before we can apply \refLemma{lemma:KushnerYin2003__Recurrence_theorem}, we need to identify a suitable function $V(\iterand{\vect{R}}{n+1})$. The choice $D(\iterand{\vect{R}}{n+1}) = \pnorm{ \iterand{\vect{R}}{n+1} - \criticalpoint{\vect{r}} }{2}^2$ comes to mind as a candidate, and we will therefore investigate \refEquation{eqn:KushnerYin2003__Inequality_to_guarantee_recurrence} for $D(\iterand{\vect{R}}{n+1})$. We will need the following result, the proof of which is relegated to \refSection{sec:Non_expansiveness_of_projection}.

\begin{lemma}
\label{lemma:Non_expansiveness_of_projection}
For $x, y \in \realNumbers$ and $\mathcal{R} = [\mathcal{R}^{\min},\mathcal{R}^{\max}] \subset \realNumbers$, $| [ x ]_{\mathcal{R}} - [ y ]_{\mathcal{R}} | \leq | x - y |$.
\end{lemma}

Combining \refEquation{eqn:Theorem_Convergence__Stochastic_gradient_algorithm} and \refLemma{lemma:Non_expansiveness_of_projection} gives
\begin{align}
D(\iterand{\vect{R}}{n+1})
&\leq \sum_{i=1}^d \bigl| \iterand{ \vectComponent{R}{i} }{n} - \itstep{n+1} \iterand{ \vectComponent{\hat{G}}{i} }{n+1} - \criticalpoint{\vectComponent{r}{i}} \bigr|^2 \nonumber \\
&= \sum_{i=1}^d \bigl| \iterand{\vectComponent{R}{i}}{n} - \criticalpoint{\vectComponent{r}{i}} \bigr|^2 + (\itstep{n+1})^2 \sum_{i=1}^d \bigl| \iterand{\vectComponent{\hat{G}}{i}}{n+1} \bigr|^2 \nonumber \\
&\phantom{=} - 2 \itstep{n+1} \sum_{i=1}^d \iterand{\vectComponent{\hat{G}}{i}}{n+1} ( \iterand{ \vectComponent{R}{i} }{n} - \criticalpoint{\vectComponent{r}{i}} ).
\end{align}
Substituting $\iterand{ \vect{\hat{G}} }{n} = \iterand{ \vect{G} }{n} + \iterand{ \vect{B} }{n} + \iterand{ \vect{E} }{n}$ into the last term, we conclude that
\begin{align}
&D(\iterand{\vect{R}}{n+1}) 
\leq D(\iterand{\vect{R}}{n}) + (\itstep{n+1})^2 \pnorm{\iterand{\vect{\hat{G}}}{n+1}}{2}^2 \nonumber \\
&- 2 \itstep{n+1} \transpose{ ( \iterand{\vect{G}}{n+1} + \iterand{\vect{B}}{n+1} + \iterand{\vect{E}}{n+1} ) } ( \iterand{\vect{R}}{n} - \criticalpoint{\vect{r}} ).
\label{eqn:Bound_on_Dn_after_nonexpansiveness}
\end{align}

Before we take the conditional expectation that results in a form similar to \refEquation{eqn:KushnerYin2003__Inequality_to_guarantee_recurrence}, recall that $u(\vect{r})$ is convex in $\vect{r}$. We therefore have that (\cite{BoydVandenberghe2004}, p.~69)
\begin{align}
&\transpose{ \iterand{\vect{G}}{n+1} } ( \criticalpoint{\vect{r}} - \iterand{\vect{R}}{n} ) 
= \transpose{ \vect{g}( \vect{\pi}(\iterand{\vect{R}}{n}), \iterand{\vect{R}}{n} ) } ( \criticalpoint{\vect{r}} - \iterand{\vect{R}}{n} ) \nonumber \\
&= \transpose{ \gradientOperatorWrt{\vect{r}} u(\iterand{\vect{R}}{n}) } ( \criticalpoint{\vect{r}} - \iterand{\vect{R}}{n} ) 
\leq u( \criticalpoint{\vect{r}} ) - u( \iterand{\vect{R}}{n} ).
\end{align} 
It follows that if $\iterand{\vect{R}}{n} \not\in \mathcal{H}_\delta$, then $\transpose{ \iterand{\vect{G}}{n+1} } ( \criticalpoint{\vect{r}} - \iterand{\vect{R}}{n} ) < - \delta / 2$. This gives in combination with \refEquation{eqn:Bound_on_Dn_after_nonexpansiveness} a term $-\delta \itstep{n+1} $, which we need for \refEquation{eqn:KushnerYin2003__Inequality_to_guarantee_recurrence}. We now note that $\expectation{ \transpose{ \iterand{\vect{E}}{n+1} } ( \iterand{\vect{R}}{n} - \criticalpoint{\vect{r}} ) | \iterand{\mathcal{F}}{n} } = 0$, so that for $\iterand{\vect{R}}{n} \not\in \mathcal{H}_\delta$, 
\begin{align}
&\expectation{ D(\iterand{\vect{R}}{n+1}) | \iterand{\mathcal{F}}{n} } - D(\iterand{\vect{R}}{n}) 
< - \delta \itstep{n+1} + \iterand{Y}{n+1}, \label{eqn:Upper_bound_for_D}
\end{align}
where 
\begin{align}
\iterand{Y}{n+1} 
= &(\itstep{n+1})^2 \expectation{  \pnorm{\iterand{\vect{\hat{G}}}{n+1}}{2}^2 | \iterand{\mathcal{F}}{n} } \nonumber \\
&+ 2 \itstep{n+1} \bigl| \expectation{ \transpose{ \iterand{\vect{B}}{n+1} } ( \iterand{\vect{R}}{n} - \criticalpoint{\vect{r}} ) | \iterand{\mathcal{F}}{n} } \bigr|.
\end{align}

The upper bound in \refEquation{eqn:Upper_bound_for_D} is not yet of the form of the right-hand side in \refEquation{eqn:KushnerYin2003__Inequality_to_guarantee_recurrence}. This implies that $D(\iterand{\vect{R}}{n+1})$ by itself is not an appropriate candidate for $V(\iterand{\vect{R}}{n+1})$. However, we can modify it slightly so that it does satisfy \refEquation{eqn:KushnerYin2003__Inequality_to_guarantee_recurrence}. For this, define $\iterand{\Delta}{n} = \expectation{ \sum_{i=n+1}^{\infty} \iterand{Y}{i} | \iterand{\mathcal{F}}{n} }$ and consider $V(\iterand{\vect{R}}{n+1}) = D(\iterand{\vect{R}}{n+1}) + \iterand{\Delta}{n+1}$ instead. The difference $\expectation{ \iterand{\Delta}{n+1} | \iterand{\mathcal{F}}{n} } - \iterand{\Delta}{n}$ is well-defined if $\sum_{i=1}^{\infty} \iterand{Y}{i} < \infty$ with probability one and is then equal to 
\begin{align}
\expectation{ \expectation{ \sum_{i=n+2}^\infty \iterand{Y}{i} | \iterand{\mathcal{F}}{n+1} } 
- \sum_{i=n+1}^\infty \iterand{Y}{i} | \iterand{\mathcal{F}}{n} } 
= - \iterand{Y}{n+1}.
\end{align}
We conclude that
\begin{align}
&\expectation{ V(\iterand{\vect{R}}{n+1}) | \iterand{\mathcal{F}}{n} } - V(\iterand{\vect{R}}{n}) \nonumber \\
&= \expectation{ D(\iterand{\vect{R}}{n+1}) | \iterand{\mathcal{F}}{n} } - D(\iterand{\vect{R}}{n}) + \expectation{ \iterand{\Delta}{n+1} | \iterand{\mathcal{F}}{n} } - \iterand{\Delta}{n} \nonumber \\
&= \expectation{ D(\iterand{\vect{R}}{n+1}) | \iterand{\mathcal{F}}{n} } - D(\iterand{\vect{R}}{n}) - \iterand{Y}{n+1}
\leq - \delta \itstep{n+1}. \label{eqn:Upper_bound_for_V}
\end{align}

The upper bound in \refEquation{eqn:Upper_bound_for_V} is of the form of \refEquation{eqn:KushnerYin2003__Inequality_to_guarantee_recurrence}, meaning that we are almost ready to apply \refLemma{lemma:KushnerYin2003__Recurrence_theorem}. What remains is to check whether 
\begin{align}
&\sum_{n=1}^{\infty} \iterand{Y}{n} 
= \sum_{n=1}^{\infty} (\itstep{n})^2 \expectation{  \pnorm{\iterand{\vect{\hat{G}}}{n}}{2}^2 | \iterand{\mathcal{F}}{n-1} } \nonumber \\
&+ 2 \sum_{n=1}^{\infty} \itstep{n} \bigl| \expectation{ \transpose{ \iterand{\vect{B}}{n} } ( \iterand{\vect{R}}{n-1} - \criticalpoint{\vect{r}} ) | \iterand{\mathcal{F}}{n-1} } \bigr| < \infty
\end{align} 
with probability one. Since $\sum_{n=1}^{\infty} (\itstep{n})^2 < \infty$ and $\pnorm{ \iterand{\vect{\hat{G}}}{n} }{2} \leq c_g$ by assumption, the first term is finite. Verifying that the second term is finite with probability one is much harder because it involves regularity conditions on $\vect{g}(\vect{\pi}(\vect{r}), \vect{r})$ and finiteness of mixing times. This can in fact be shown as stated in the next lemma, proved in \refSection{sec:Error_bias}.

\begin{lemma}
\label{lemma:Recurrence_of_A_when_summing_Yi_is_finite_wp1}
Under the assumptions of \refTheorem{thm:Convergence}, the sum $\sum_{n=1}^{\infty} \itstep{n} \bigl| \expectation{ \transpose{ \iterand{\vect{B}}{n} } ( \iterand{\vect{R}}{n-1} - \criticalpoint{\vect{r}} ) | \iterand{\mathcal{F}}{n-1} } \bigr|$ is finite with probability one.
\end{lemma}

\subsubsection{Capture}

Having derived conditions under which $\mathcal{H}_{\delta}$ is recurrent, we turn our attention to deriving conditions under which capture occurs. Recall that capture means that there must exist an $m \in \naturalNumbers$ large enough so that $D( \iterand{\vect{R}}{n} ) \leq D(\iterand{\vect{R}}{m}) + \varepsilon$ for all $n \geq m$ with probability one. 

After applying \refEquation{eqn:Bound_on_Dn_after_nonexpansiveness} repeatedly and using the upper bound $\transpose{ \iterand{\vect{G}}{n} } ( \iterand{\vect{R}}{n-1} - \criticalpoint{\vect{r}} ) \geq 0$, which follows by convexity of $u(\vect{r})$, we find that
\begin{align}
&D(\iterand{\vect{R}}{n}) 
\leq D(\iterand{\vect{R}}{m}) + \sum_{j=m}^{n} (\itstep{j+1})^2 \pnorm{\iterand{\vect{\hat{G}}}{j+1}}{2}^2
\nonumber \\
&- 2 \sum_{j=m}^{n} \itstep{j+1} \transpose{ ( \iterand{\vect{B}}{j+1} + \iterand{\vect{E}}{j+1} ) } ( \iterand{\vect{R}}{j} - \criticalpoint{\vect{r}} ).
\label{eqn:Repeated_application_of_the_bound_on_Dn_after_nonexpansiveness}
\end{align}
We now need to show that each sum in the right-hand side of \refEquation{eqn:Repeated_application_of_the_bound_on_Dn_after_nonexpansiveness} becomes small for $m$ sufficiently large. Because $\sum_{n=1}^{\infty} (\itstep{n})^2 < \infty$ and $\pnorm{ \iterand{\vect{\hat{G}}}{n} }{2} \leq c_g$, it immediately follows that $\lim_{m \rightarrow \infty} \sum_{j=m}^{\infty} (\itstep{n})^2 \pnorm{ \iterand{\vect{\hat{G}}}{n} }{2} = 0$. In turn, this implies that for any $\varepsilon$, there exists an $m_{0} \in \naturalNumbers$ so that $\sum_{j=m}^{n} (\itstep{n})^2 \pnorm{ \iterand{\vect{\hat{G}}}{n} }{2} \leq \varepsilon$ for all $n \geq m \geq m_{0}$. Verifying that the other two sums become small is substantially more difficult. This can be established using martingale arguments, as asserted in \refLemma{lemma:Capture_of_R_when_close_enough}, the proof of which is postponed to \refSection{sec:Zero_mean_noise}.

\begin{lemma}
\label{lemma:Capture_of_R_when_close_enough}
Under the assumptions of \refTheorem{thm:Convergence}, for any $\varepsilon > 0$, there exists $m_{0} \in \naturalNumbers$ so that for any $n \geq m \geq m_{0}$

\textnormal{(i)} $\sum_{j=m}^{n} \itstep{j} \transpose{ \iterand{\vect{B}}{j} } ( \criticalpoint{\vect{r}} - \iterand{\vect{R}}{j-1} ) \leq \varepsilon$ and

\textnormal{(ii)} $\sum_{j=m}^{n} \itstep{j} \transpose{ \iterand{\vect{E}}{j} } ( \criticalpoint{\vect{r}} - \iterand{\vect{R}}{j-1} ) \leq \varepsilon$

\noindent with probability one.
\end{lemma}

Our work thus far can also be used to prove that the gradient algorithm \refEquation{eqn:Analytic_gradient_algorithm} converges. It is a special case of its stochastic counterpart \refEquation{eqn:Theorem_Convergence__Stochastic_gradient_algorithm}, for which $\iterand{\vect{B}}{n} = \vect{0}$, $\iterand{\vect{E}}{n} = \vect{0}$, $\iterand{\vect{G}}{n} = \iterand{\vect{\hat{G}}}{n} = \iterand{\vect{g}}{n}$ and $\iterand{\vect{R}}{n} = \iterand{\vect{r}}{n}$ for all $n \geq 0$. To prove that \refEquation{eqn:Analytic_gradient_algorithm} converges, we apply \refEquation{eqn:Bound_on_Dn_after_nonexpansiveness} repeatedly and use that $\transpose{ \iterand{\vect{g}}{n} } ( \criticalpoint{\vect{r}} - \iterand{\vect{r}}{n-1} ) \leq u( \criticalpoint{\vect{r}} ) - u( \iterand{\vect{r}}{n-1} )$ for any $n \in \naturalNumbers$ by convexity of $u(\vect{r})$, so that
\begin{align}
D(\iterand{\vect{r}}{n}) 
&\leq D(\iterand{\vect{r}}{0}) + \sum_{j=0}^{n} (\itstep{j+1})^2 \pnorm{\iterand{\vect{g}}{j+1}}{2}^2 \nonumber \\
&\phantom{=} - 2 \sum_{j=0}^{n} \itstep{j+1} ( u( \iterand{\vect{r}}{j} ) - u( \criticalpoint{\vect{r}} ) ).
\end{align}
Noting that $D(\iterand{\vect{r}}{n}) \geq 0$ for all $n \in \naturalNumbers$ and $D(\iterand{\vect{r}}{0}) \leq \startingErrorConstant$ for some constant $\startingErrorConstant < \infty$ since $\iterand{\vect{r}}{0} \in \mathcal{R}$, we conclude that
\begin{align}
2 \sum_{j=0}^{n} \itstep{j+1} ( u( \iterand{\vect{r}}{j} ) - u( \criticalpoint{\vect{r}} ) ) \leq \startingErrorConstant + \gradientConstant^2  \sum_{j=0}^{n} (\itstep{j+1})^2.
\end{align}
Since $\sum_{j=0}^{n} \itstep{j+1} ( u( \iterand{\vect{r}}{j} ) - u( \criticalpoint{\vect{r}} ) ) \geq \min_{i=0,...,n} \{ u( \iterand{\vect{r}}{i} ) - u( \criticalpoint{\vect{r}} ) \} \sum_{j=0}^{n} \itstep{j+1}$, we have the inequality
\begin{align}
\min_{i=0,...,n} \{ u( \iterand{\vect{r}}{i} ) - u( \criticalpoint{\vect{r}} ) \}
\leq \frac{ \startingErrorConstant + \gradientConstant^2  \sum_{j=0}^{n} (\itstep{j+1})^2 }{ 2 \sum_{j=0}^{n} \itstep{j+1} },
\end{align}
which converges to $0$ as $n \rightarrow \infty$. 

From this little detour we see that it is much easier to establish convergence for \refEquation{eqn:Analytic_gradient_algorithm} than for its stochastic counterpart \refEquation{eqn:Theorem_Convergence__Stochastic_gradient_algorithm}. It is the error bias and zero-mean noise that make the convergence analysis of \refEquation{eqn:Theorem_Convergence__Stochastic_gradient_algorithm} so much harder.
 
\subsection{Evaluating the conditions}
\label{sec:Proving_convergence}

We now provide the proofs of \refLemma{lemma:Recurrence_of_A_when_summing_Yi_is_finite_wp1} and \ref{lemma:Capture_of_R_when_close_enough}, which together prove \refTheorem{thm:Convergence}. In our proofs, we choose to consider the error bias and zero-mean noise separately, which makes the analysis more tractable.

\subsubsection{Error bias}
\label{sec:Error_bias}

We start by showing that the error bias satisfies the property claimed in \refLemma{lemma:Recurrence_of_A_when_summing_Yi_is_finite_wp1} under the assumptions of \refTheorem{thm:Convergence}. After substituting the definition of the error bias and using the triangle inequality, one finds that
\begin{align}
&\sum_{n=1}^{\infty} \itstep{n} \bigl| \expectation{ \transpose{ \iterand{\vect{B}}{n} } ( \iterand{\vect{R}}{n-1} - \criticalpoint{\vect{r}} ) | \iterand{\mathcal{F}}{n-1} } \bigr| \nonumber \\
&= \sum_{n=1}^{\infty} \itstep{n} \bigl| \sum_{i=1}^d \expectation{ \iterand{\vectComponent{B}{i}}{n} | \iterand{\mathcal{F}}{n-1} } ( \iterand{\vectComponent{R}{i}}{n-1} - \criticalpoint{\vectComponent{r}{i}} ) \bigr| \nonumber \\
&\leq \sum_{n=1}^{\infty} \itstep{n} \sum_{i=1}^d ( \vectComponent{\mathcal{R}^{\max}}{i} - \vectComponent{\mathcal{R}^{\min}}{i} ) \bigl| \expectation{ \iterand{\vectComponent{B}{i}}{n} | \iterand{\mathcal{F}}{n-1} } \bigr| \nonumber \\
&= \sum_{n=1}^{\infty} \itstep{n} \sum_{i=1}^d ( \vectComponent{\mathcal{R}^{\max}}{i} - \vectComponent{\mathcal{R}^{\min}}{i} ) \bigl| \iterand{\vectComponent{B}{i}}{n} \bigr|.
\label{eqn:Bound_on_the_expectation_of_B_times_R_minus_rstar__up_to_abs_of_B}
\end{align}
The inequality is a consequence of $\mathcal{R}$ being a hypercube. We have also used the fact that $\expectation{ \iterand{\vectComponent{B}{i}}{n} | \iterand{\mathcal{F}}{n-1} } = \iterand{\vectComponent{B}{i}}{n}$, which follows from the definition $\iterand{\vectComponent{G}{i}}{n} = \vectComponent{g}{i}( \vect{\pi}(\iterand{\vect{R}}{n-1}), \iterand{\vect{R}}{n-1} )$.

We now bound $\bigl| \iterand{\vectComponent{B}{i}}{n} \bigr|$ from above. After recalling that $\iterand{\vectComponent{B}{i}}{n} = \expectation{ \iterand{ \vectComponent{\hat{G}}{i} }{n} | \iterand{\mathcal{F}}{n-1} } - \iterand{ \vectComponent{G}{i} }{n}$ and using Jensen's inequality, we find that $\bigl| \iterand{\vectComponent{B}{i}}{n} \bigr|$ equals
\begin{align}
&\bigl| \expectation{ \vectComponent{g}{i}( \iterand{\vect{\hat{\Pi}}}{n}, \iterand{\vect{R}}{n-1} ) | \iterand{\mathcal{F}}{n-1} } - \vectComponent{g}{i}( \vect{\pi}(\iterand{\vect{R}}{n-1}), \iterand{\vect{R}}{n-1} ) \bigr| \nonumber \\
&= \bigl| \expectation{ \vectComponent{g}{i}( \iterand{\vect{\hat{\Pi}}}{n}, \iterand{\vect{R}}{n-1} ) - \vectComponent{g}{i}( \vect{\pi}(\iterand{\vect{R}}{n-1}), \iterand{\vect{R}}{n-1} ) | \iterand{\mathcal{F}}{n-1} } \bigr| \nonumber \\
&\leq \expectation{ \bigl| \vectComponent{g}{i}( \iterand{\vect{\hat{\Pi}}}{n}, \iterand{\vect{R}}{n-1} ) - \vectComponent{g}{i}( \vect{\pi}(\iterand{\vect{R}}{n-1}), \iterand{\vect{R}}{n-1} ) \bigr| | \iterand{\mathcal{F}}{n-1} }. \nonumber
\end{align}
Recalling condition \refEquation{eqn:Theorem_Convergence__Assumption_Lipschitz_continuity} gives
\begin{align}
\bigl| \iterand{\vectComponent{B}{i}}{n} \bigr|
\leq \frac{\LipschitzConstant}{2} \sum_{\svA \in \stateSpace} \expectation{ \bigl| \iterand{\vectComponent{\hat{\Pi}}{\svA}}{n} - \vectComponent{\pi}{\svA}(\iterand{\vect{R}}{n-1}) \bigr| | \iterand{\mathcal{F}}{n-1} }. \label{eqn:Error_bias__Upper_bound_after_triangle_inequality}
\end{align}

Finiteness of \refEquation{eqn:Bound_on_the_expectation_of_B_times_R_minus_rstar__up_to_abs_of_B} can now be proven by constructing an upper bound for \refEquation{eqn:Error_bias__Upper_bound_after_triangle_inequality}. We can obtain such a bound using the following lemma, proved in \refAppendixSection{sec:Appendix__Proving_the_probabilistic_error_bound_on_the_steady_state_estimate}.

\begin{lemma}
\label{lemma:Ergodic_probability_bound}
There exist $\ergodicConstant, \kappa \in [0,\infty)$ such that for $\iterror{n} \in [0,1]$ and $\svA \in \stateSpace$,
\begin{align}
\probability{ \bigl| \iterand{\vectComponent{\hat{\Pi}}{\svA}}{n} - \vectComponent{\pi}{\svA}(\iterand{\vect{R}}{n-1}) \bigr| \geq \iterror{n} } 
\leq \ergodicConstant \exp{ \Bigl( - \frac{ (\iterror{n})^2 }{ 4 \cardinality{\stateSpace}^2 \kappa \itfreq{n} } \Bigr) }. \nonumber
\end{align}
\end{lemma}

Define $\iterand{\vectComponent{\Phi}{\svA}}{n} = \bigl| \iterand{\vectComponent{\hat{\Pi}}{\svA}}{n} - \vectComponent{\pi}{\svA}(\iterand{\vect{R}}{n-1}) \bigr|$ and let $\iterand{\epsilon}{n} \in [0,1]$. Using \refEquation{eqn:Error_bias__Upper_bound_after_triangle_inequality} and then \refLemma{lemma:Ergodic_probability_bound} yields
\begin{align}
\bigl| \iterand{\vectComponent{B}{i}}{n} \bigr|
&\leq \frac{\LipschitzConstant}{2} \sum_{\svA \in \stateSpace} \expectation{ \iterand{\vectComponent{\Phi}{\svA}}{n} | \iterand{\mathcal{F}}{n-1} } \nonumber \\
&= \frac{\LipschitzConstant}{2} \sum_{\svA \in \stateSpace} \Bigl( 
\probability{ \iterand{\vectComponent{\Phi}{\svA}}{n} < \iterror{n} } \expectation{ \iterand{\vectComponent{\Phi}{\svA}}{n} | \iterand{\mathcal{F}}{n-1}, \iterand{\vectComponent{\Phi}{\svA}}{n} < \iterror{n} } \nonumber \\
&\phantom{=} + \probability{ \iterand{\vectComponent{\Phi}{\svA}}{n} \geq \iterror{n} } \expectation{ \iterand{\vectComponent{\Phi}{\svA}}{n} | \iterand{\mathcal{F}}{n-1}, \iterand{\vectComponent{\Phi}{\svA}}{n} \geq \iterror{n} } \Bigr) \nonumber \\
&\leq \frac{\LipschitzConstant}{2} \sum_{\svA \in \stateSpace} \Bigl( 
\iterror{n} + ( 1 - \iterror{n} ) \probability{ \iterand{\vectComponent{\Phi}{\svA}}{n} \geq \iterror{n} } \Bigr) \nonumber \\
&\leq \frac{\LipschitzConstant \cardinality{\stateSpace}}{2} \max \{ 1, \ergodicConstant \} \Bigl( \iterror{n} + \exp{ \Bigl( - \frac{ (\iterror{n})^2 }{ 4 \cardinality{\stateSpace}^2 \kappa \itfreq{n} } \Bigr) } \Bigr). \label{eqn:Upper_bound_for_Bn_after_probability_estimate}
\end{align}
After bounding \refEquation{eqn:Bound_on_the_expectation_of_B_times_R_minus_rstar__up_to_abs_of_B} from above using \refEquation{eqn:Upper_bound_for_Bn_after_probability_estimate}, it follows from \refEquation{eqn:Theorem_Convergence__Assumption_decreasing_error_and_update_frequency} that 
\begin{align}
\sum_{n=1}^{\infty} \itstep{n} \bigl| \expectation{ \transpose{ \iterand{\vect{B}}{n} } ( \iterand{\vect{R}}{n-1} - \criticalpoint{\vect{r}} ) | \iterand{\mathcal{F}}{n-1} } \bigr| < \infty,
\end{align}
which completes the proof of \refLemma{lemma:Recurrence_of_A_when_summing_Yi_is_finite_wp1}.

We now show that the error bias satisfies assertion (i) in \refLemma{lemma:Capture_of_R_when_close_enough} under the assumptions of \refTheorem{thm:Convergence}. Similar to the derivation of \refEquation{eqn:Bound_on_the_expectation_of_B_times_R_minus_rstar__up_to_abs_of_B},
\begin{align}
&\sum_{n=1}^{\infty} \itstep{n} \bigl| \transpose{ \iterand{\vect{B}}{n} } ( \iterand{\vect{R}}{n-1} - \criticalpoint{\vect{r}} ) \bigr| \nonumber \\
&\leq \sum_{n=1}^{\infty} \itstep{n} \sum_{i=1}^d ( \vectComponent{\mathcal{R}^{\max}}{i} - \vectComponent{\mathcal{R}^{\min}}{i} ) \bigl| \iterand{\vectComponent{B}{i}}{n} \bigr|. \label{eqn:Error_bias__Intermediate_step_A}
\end{align}
Combining \refEquation{eqn:Error_bias__Intermediate_step_A}, \refEquation{eqn:Upper_bound_for_Bn_after_probability_estimate} and \refEquation{eqn:Theorem_Convergence__Assumption_decreasing_error_and_update_frequency}, we conclude that with probability one, 
\begin{align}
\sum_{n=1}^{\infty} \itstep{n} \bigl| \transpose{ \iterand{\vect{B}}{n} } ( \iterand{\vect{R}}{n-1} - \criticalpoint{\vect{r}} ) \bigr| < \infty,
\end{align}
so that $\lim_{m \rightarrow \infty} \sum_{j=m}^{\infty} \itstep{n} \transpose{ \iterand{\vect{B}}{n} } ( \criticalpoint{\vect{r}} - \iterand{\vect{R}}{n-1} ) = 0$ with probability one. This implies that there exists an $m_{0} \in \naturalNumbers$ so that for all $n \geq m \geq m_{0}$, $\sum_{j=m}^{n} \itstep{j} \transpose{ \iterand{\vect{B}}{j} } ( \criticalpoint{\vect{r}} - \iterand{\vect{R}}{j-1} ) \leq \varepsilon$ with probability one. The error bias thus satisfies assertion (i) in \refLemma{lemma:Capture_of_R_when_close_enough}. All that remains is to show that the zero-mean noise satisfies \refLemma{lemma:Capture_of_R_when_close_enough}(ii).

\subsubsection{Zero-mean noise}
\label{sec:Zero_mean_noise}

We use a martingale argument to show that assertion (ii) in \refLemma{lemma:Capture_of_R_when_close_enough} holds. We start our argument by defining $\iterand{M}{n} = \sum_{j=1}^n \iterand{a}{j} \transpose{ \iterand{\vect{E}}{j} } ( \iterand{\vect{R}}{j-1} - \criticalpoint{\vect{r}} )$. See \refAppendixSection{sec:Appendix__Proof_that_Mn_is_a_martingale} for a proof of the following result.

\begin{lemma}
\label{lemma:Mn_is_a_martingale}
$\iterand{M}{n}$ is a martingale.
\end{lemma}

We will use a martingale convergence theorem \cite{Steele2001} to show that for $n \geq m$ both sufficiently large, $\iterand{M}{n} - \iterand{M}{m-1} \leq \varepsilon$ with probability one.

\begin{theorem}
\label{thm:Martingale_convergence}
If $\{ \iterand{M}{n} \}$ is a martingale for which there exists a constant $\martingaleConstant < \infty$ so that $\expectation{ (\iterand{M}{n})^2 } \leq \martingaleConstant$ for all $n \geq 0$, then there exists a random variable $\criticalpoint{M}$ with $\expectation{ (\criticalpoint{M})^2 } \leq \martingaleConstant$ such that $\iterand{M}{n} \rightarrow \criticalpoint{M}$ with probability one as $n \rightarrow \infty$. Moreover, $\expectation{ | \iterand{M}{n} - \criticalpoint{M} |^2 }^{\frac{1}{2}} \rightarrow 0$ as $n \rightarrow \infty$.
\end{theorem}

Before we can apply \refTheorem{thm:Martingale_convergence}, we need to show existence of a $\martingaleConstant \in \realNumbers$ such that $\expectation{ (\iterand{M}{n})^2 } \leq \martingaleConstant$ for all $n \in \naturalNumbers$. To show this, expand
\begin{align}
&\sup_{n} \expectation{ (\iterand{M}{n})^2 }
= \sup_{n} \bigl\{ \sum_{j=1}^n (\itstep{j})^2 \expectation{ ( \transpose{ \iterand{\vect{E}}{j} } ( \iterand{\vect{R}}{j-1} - \criticalpoint{\vect{r}} ) )^2 } \nonumber \\
&+ \sum_{ j \neq k } \iterand{a}{j} \iterand{a}{k} \expectation{ \transpose{ \iterand{\vect{E}}{j} } ( \iterand{\vect{R}}{j-1} - \criticalpoint{\vect{r}} ) \transpose{ \iterand{\vect{E}}{k} } ( \iterand{\vect{R}}{k-1} - \criticalpoint{\vect{r}} ) } \bigr\}, \nonumber
\end{align}
and then consider any one of the cross terms with $k < j$. By the tower property,
\begin{align}
&\, \expectation{ \transpose{ \iterand{\vect{E}}{j} } ( \iterand{\vect{R}}{j-1} - \criticalpoint{\vect{r}} ) \transpose{ \iterand{\vect{E}}{k} } ( \iterand{\vect{R}}{k-1} - \criticalpoint{\vect{r}} ) } \nonumber \\
= &\, \expectation{ \expectation{ \transpose{ \iterand{\vect{E}}{j} } ( \iterand{\vect{R}}{j-1} - \criticalpoint{\vect{r}} ) \transpose{ \iterand{\vect{E}}{k} } ( \iterand{\vect{R}}{k-1} - \criticalpoint{\vect{r}} ) | \iterand{\mathcal{F}}{j-1} } } \nonumber \\
= &\, \expectation{ \transpose{ \iterand{\vect{E}}{k} } ( \iterand{\vect{R}}{k-1} - \criticalpoint{\vect{r}} ) \expectation{ \sum_{i=1}^d \iterand{\vectComponent{E}{i}}{j} ( \iterand{\vectComponent{R}{i}}{j-1} - \criticalpoint{\vectComponent{r}{i}} ) | \iterand{\mathcal{F}}{j-1} } } \nonumber \\
= &\, \expectation{ \transpose{ \iterand{\vect{E}}{k} } ( \iterand{\vect{R}}{k-1} - \criticalpoint{\vect{r}} ) \sum_{i=1}^d \expectation{ \iterand{\vectComponent{E}{i}}{j} | \iterand{\mathcal{F}}{j-1} } ( \iterand{\vectComponent{R}{i}}{j-1} - \criticalpoint{\vectComponent{r}{i}} ) }, \nonumber
\end{align}
and because $\expectation{ \iterand{\vectComponent{E}{i}}{j} | \iterand{\mathcal{F}}{j-1} } = 0$, all cross terms are equal to $0$. Because the summands are positive, we can give an upper bound by summing over all terms, so that
\begin{align}
&\sup_{n} \expectation{ (\iterand{M}{n})^2 } 
\leq \sum_{j=1}^\infty (\itstep{j})^2 \expectation{ ( \transpose{ \iterand{\vect{E}}{j} } ( \iterand{\vect{R}}{j-1} - \criticalpoint{\vect{r}} ) )^2 } \nonumber \\
&= \sum_{j=1}^\infty (\itstep{j})^2 \expectation{ \bigl( \sum_{i=1}^d \iterand{\vectComponent{E}{i}}{j} ( \iterand{\vectComponent{R}{i}}{j-1} - \criticalpoint{\vectComponent{r}{i}} ) \bigr)^2 }.
\end{align}
Using the triangle inequality, we find that
\begin{align}
\sup_{n} \expectation{ (\iterand{M}{n})^2 }
\leq &\sum_{j=1}^\infty (\itstep{j})^2 \expectation{ \bigl( \sum_{i=1}^d | \iterand{\vectComponent{E}{i}}{j} | | \iterand{\vectComponent{R}{i}}{j-1} - \criticalpoint{\vectComponent{r}{i}} | \bigr)^2 }. \nonumber
\end{align}
Now note that $\sum_{i=1}^d \bigl| \iterand{\vectComponent{E}{i}}{j} \bigr| = \pnorm{ \iterand{ \vect{E} }{j} }{1}$, write
\begin{align}
&\pnorm{ \iterand{ \vect{E} }{j} }{1}
\leq \expectation{ \pnorm{ \iterand{ \vect{\hat{G}} }{j} }{1} | \iterand{\mathcal{F}}{j-1} }
+ \pnorm{ \iterand{ \vect{\hat{G}} }{j} }{1} \nonumber \\
&\leq \sqrt{d} \expectation{ \pnorm{ \iterand{ \vect{\hat{G}} }{j} }{2} | \iterand{\mathcal{F}}{j-1} } + \sqrt{d} \pnorm{ \iterand{ \vect{\hat{G}} }{j} }{2} 
\leq 2 \gradientConstant \sqrt{d} \label{eqn:Error_bias__Upper_bound_in_terms_of_dimension_and_gradient}
\end{align}
and recall that $\mathcal{R}$ is a hypercube. We conclude that
\begin{align}
\sup_{n} \expectation{ (\iterand{M}{n})^2 }
\leq 4 \gradientConstant^2 d \max_{i=1,...,d} \{ ( \vectComponent{\mathcal{R}^{\max}}{i} - \vectComponent{\mathcal{R}^{\min}}{i} )^2 \} \sum_{j=1}^\infty (\itstep{j})^2. \nonumber
\end{align}
The right-hand side is finite by condition \refEquation{eqn:Theorem_Convergence__Assumption_decreasing_stepsizes}, and we see that there indeed exists a coefficient $\martingaleConstant$ so that $\expectation{ (\iterand{M}{n})^2 } \leq \martingaleConstant$ for all $n \in \naturalNumbers$. We now apply \refTheorem{thm:Martingale_convergence} and conclude that as $n \geq m \rightarrow \infty$,
\begin{align}
&\expectation{ | \iterand{M}{n} - \iterand{M}{m-1} |^2 }^{\frac{1}{2}}
\leq \expectation{ | \iterand{M}{n} - \criticalpoint{M} |^2 }^{\frac{1}{2}} \nonumber \\
&+ \expectation{ | \iterand{M}{m-1} - \criticalpoint{M} |^2 }^{\frac{1}{2}} \rightarrow 0.
\end{align}

This result enables us to use Doob's maximal inequality \cite{Steele2001}, as reproduced in the lemma below, in order to conclude that \refLemma{lemma:Capture_of_R_when_close_enough}(ii) holds.

\begin{lemma}
\label{lemma:Doobs_maximal_inequality}
If $\process{ \iterand{M}{n} }{ n \geq 0 }$ is a nonnegative submartingale and $\lambda > 0$, then
\begin{align}
\lambda \probability{ \sup_{m \leq n} \iterand{M}{m} \geq \lambda }
\leq \expectation{ \iterand{M}{n} \indicator{ \sup_{m \leq n} \iterand{M}{m} \geq \lambda } }
\leq \expectation{ \iterand{M}{n} }. \nonumber
\end{align}
\end{lemma}

Fix $m \in \naturalNumbers$ and define $\iterand{W}{n} = \iterand{M}{n+m-1} - \iterand{M}{m-1}$ for $n \in \naturalNumbers$. $|\iterand{W}{n}|$ is a submartingale by Jensen's inequality with respect to the sequence $\iterand{\mathcal{F}}{m-1}, \iterand{\mathcal{F}}{m}, \iterand{\mathcal{F}}{m+1}, ...~$, since $\expectation{ | \iterand{W}{n+1} | | \iterand{\mathcal{F}}{n+m-1} } \geq \expectation{ \iterand{W}{n+1} | \iterand{\mathcal{F}}{n+m-1} } = \iterand{W}{n}$. Applying \refLemma{lemma:Doobs_maximal_inequality} to $|\iterand{W}{n}|$, we find that
\begin{align}
&\probability{ \sup_{0 \leq t \leq n} | \iterand{M}{t+m-1} - \iterand{M}{m-1} | \geq \lambda } \nonumber \\
&\leq \frac{ \expectation{ | \iterand{M}{n+m-1} - \iterand{M}{m-1} | } }{ \lambda } \nonumber \\
&\leq \frac{ \expectation{ | \iterand{M}{n+m-1} - \criticalpoint{M} | } + \expectation{ | \criticalpoint{M} - \iterand{M}{m-1} | } }{ \lambda } \nonumber \\
&\leq \frac{ \expectation{ | \iterand{M}{n+m-1} - \criticalpoint{M} |^2 }^{\frac{1}{2}} + \expectation{ | \criticalpoint{M} - \iterand{M}{m-1} |^2 }^{\frac{1}{2}} }{ \lambda } \nonumber
\end{align}
for any $\lambda \in \strictlyPositiveRealNumbers$ and $m \in \naturalNumbers$. This upper bound converges to $0$ as $n, m \rightarrow \infty$, implying that there exists an $m_{0} \in \naturalNumbers$, such that for all $n \geq m \geq m_{0}$, $\iterand{M}{n} - \iterand{M}{m-1} \leq \varepsilon$ with probability one.

Having established \refLemma{lemma:Recurrence_of_A_when_summing_Yi_is_finite_wp1} and \ref{lemma:Capture_of_R_when_close_enough}, the proof of \refTheorem{thm:Convergence} is now completed.

\section{Conclusions}

We have developed an online gradient algorithm for finding parameter values that optimize the performance of reversible Markov processes with product-form distributions. As a key feature, the approach avoids the computational complexity of calculating the gradient in terms of the stationary probabilities and instead relies on measuring empirical time fractions of the various states so as to obtain estimates for the gradient. While the impact of the induced measurement noise can be handled without too much trouble, the bias in the estimates presents a trickier issue. In order to exploit mixing time results to deal with the bias, we focussed on reversible processes. We expect however that convergence can be established under milder conditions.

For fast convergence, the algorithm needs to strike a balance between the step sizes and the lengths of observation periods, which is a consequence of the existence of two time scales - one being the mixing time of the underlying stochastic process and the other being the iteration sequence generated by the algorithm. Intuitively, the step sizes should not have become too small by the time that the observation periods have become larger than the mixing time. The convergence of the algorithm would otherwise slow down drastically. A challenging issue for further research is to gain a more detailed understanding of the effect of step sizes and the role of mixing times in relation to the convergence speed. A related direction is to explore the trade-off between accuracy in static scenarios and responsiveness in dynamic environments, which relates to convergence in distribution for non-vanishing step sizes as opposed to the almost-sure convergence for decreasing step sizes as considered here.

\section*{Acknowledgments}

This research was financially supported by The Netherlands Organization for Scientific Research (NWO) in the framework of the TOP-GO program and by an ERC Starting Grant.

\bibliographystyle{IEEEtran}
\bibliography{IEEEabrv,Bibliography}

\appendix

\section{Remaining proofs}

\subsection{Proof of Lemma \ref{lemma:Lipschitz_continuity_of_g_for_MMss_queue}}
\label{sec:Appendix__Lipschitz_continuity}

Define $\ErlangBwrt{\vect{\mu}} = \vectComponent{\mu}{\numberOfServers}$ and $\ExpectedBwrt{\vect{\mu}} = \sum_{\sA = 1}^{\numberOfServers} \sA \vectComponent{\mu}{\sA}$ for all $\vect{\mu} \in [0,1]^{\cardinality{\stateSpace}}$ for which $\transpose{ \vectones{\cardinality{\stateSpace}} } \vect{\mu} = 1$. By definition of $g(\vect{\mu},r)$, $\pnorm{ g(\vect{\mu},r) }{2} \leq | \ErlangBwrt{\vect{\mu}} ( \ExpectedBwrt{\vect{\mu}} - \numberOfServers ) | / r + | c'(r) | < \infty$.
The first term is finite because $r \geq \mathcal{R}^{\min} > 0$, $\ErlangBwrt{\vect{\mu}} \leq 1$ and $\ExpectedBwrt{\vect{\mu}} \leq s < \infty$. The second term is finite by our assumption that $c'(r)$ is bounded for all $r \in \mathcal{R}$. This proves that condition \refEquation{eqn:Theorem_Convergence__Assumption_Boundedness_of_g_in_r} is met. 

We now turn to condition \refEquation{eqn:Theorem_Convergence__Assumption_Lipschitz_continuity}. Write
$| g( \vect{\mu}, r ) - g( \vect{\nu}, r ) | = | \ErlangBwrt{\vect{\mu}} ( \ExpectedBwrt{\vect{\mu}} - \numberOfServers ) - \ErlangBwrt{\vect{\nu}} ( \ExpectedBwrt{\vect{\nu}} - \numberOfServers ) | / r \leq | \ErlangBwrt{\vect{\mu}} \ExpectedBwrt{\vect{\mu}} - \numberOfServers \ErlangBwrt{\vect{\mu}} - \ErlangBwrt{\vect{\nu}} \ExpectedBwrt{\vect{\nu}} + \numberOfServers \ErlangBwrt{\vect{\nu}} | / \mathcal{R}^{\min} \leq ( | \ErlangBwrt{\vect{\mu}} \ExpectedBwrt{\vect{\mu}} - \ErlangBwrt{\vect{\nu}} \ExpectedBwrt{\vect{\nu}} | + \numberOfServers | \ErlangBwrt{\vect{\mu}} - \ErlangBwrt{\vect{\nu}} | ) / \mathcal{R}^{\min}$. We then conclude that $| \ErlangBwrt{\vect{\mu}} \ExpectedBwrt{\vect{\mu}} - \ErlangBwrt{\vect{\nu}} \ExpectedBwrt{\vect{\nu}} | = | \ErlangBwrt{\vect{\mu}} \ExpectedBwrt{\vect{\mu}} -  \ErlangBwrt{\vect{\mu}} \ExpectedBwrt{\vect{\nu}} + \ErlangBwrt{\vect{\mu}} \ExpectedBwrt{\vect{\nu}} - \ErlangBwrt{\vect{\nu}} \ExpectedBwrt{\vect{\nu}} | \leq \ErlangBwrt{\vect{\mu}} | \ExpectedBwrt{\vect{\mu}} - \ExpectedBwrt{\vect{\nu}} | + \ExpectedBwrt{\vect{\nu}} | \ErlangBwrt{\vect{\mu}} - \ErlangBwrt{\vect{\nu}} | \leq | \ExpectedBwrt{\vect{\mu}} - \ExpectedBwrt{\vect{\nu}} | + \numberOfServers | \ErlangBwrt{\vect{\mu}} - \ErlangBwrt{\vect{\nu}} |$,
so that $| g( \vect{\mu}, r ) - g( \vect{\nu}, r ) | 
\leq ( | \ExpectedBwrt{\vect{\mu}} - \ExpectedBwrt{\vect{\nu}} | + 2 \numberOfServers | \ErlangBwrt{\vect{\mu}} - \ErlangBwrt{\vect{\nu}} | ) / \mathcal{R}^{\min}$. Finally, by definition of $\ErlangBwrt{\vect{\mu}}$, $
| \ErlangBwrt{\vect{\mu}} - \ErlangBwrt{\vect{\nu}} | = | \vectComponent{\mu}{\numberOfServers} - \vectComponent{\nu}{\numberOfServers} | \leq 2 \totalVariation{ \vect{\mu} - \vect{\nu} }$.
Similarly for $\ExpectedBwrt{\vect{\mu}}$,
$\bigl| \ExpectedBwrt{\vect{\mu}} - \ExpectedBwrt{\vect{\nu}} \bigr| 
\leq \sum_{\sA=1}^{\numberOfServers} \sA \bigl| \vectComponent{\mu}{\sA} -  \vectComponent{\nu}{\sA} \bigr|
\leq 2 \numberOfServers \totalVariation{ \vect{\mu} - \vect{\nu} }$. Thus $| g( \vect{\mu}, r ) - g( \vect{\nu}, r ) | \leq 6 \numberOfServers \totalVariation{ \vect{\mu} - \vect{\nu} } / \mathcal{R}^{\min}$,
which concludes the proof after setting $\LipschitzConstant = 6 \numberOfServers / \mathcal{R}^{\min}$. \QuodEratDemonstrandum

\subsection{Proof of Lemma \ref{lemma:Convexity_of_the_log_likelihood_function}}
\label{sec:Appendix__Convexity_of_the_log_likelihood_function}

Substituting \refEquation{eqn:Product_form__Equilibrium_distribution} into \refEquation{eqn:Introduction_loglikelihood_function} gives
\begin{align}
u(\vect{r}) 
= \ln \sum_{\svB \in \stateSpace} \exp{ \vectComponent{ ( A \vect{r} + \vect{b} ) }{\svB} } - \sum_{\svA \in \stateSpace} \vectComponent{\alpha}{\svA} \vectComponent{ ( A \vect{r} + \vect{b} ) }{\svA}.
\end{align}
The function $v(\vect{s}) = \ln \sum_{\svB \in \stateSpace} \exp{ \vectComponent{s}{\svB} } - \sum_{\svA \in \stateSpace} \vectComponent{\alpha}{\svA} \vectComponent{s}{\svA}$ is convex on $\realNumbers^{\cardinality{\stateSpace}}$ \cite{BoydVandenberghe2004}, p.~72. We see that $u(\vect{r})$ is a composition of a convex function with an affine mapping, i.e.~$u(\vect{r}) = v(A\vect{r}+\vect{b})$, and such functions are convex \cite{BoydVandenberghe2004}, p.~79. \QuodEratDemonstrandum

\subsection{Proof of Lemma \ref{lemma:Non_expansiveness_of_projection}}
\label{sec:Non_expansiveness_of_projection}

Define $l = \mathcal{R}^{\min}$ and $r = \mathcal{R}^{\max}$. If $x, y \in \mathcal{R}$, equality holds. Consider the case $x \not\in \mathcal{R}, y \in \mathcal{R}$. If $x > r$, $| [ x ]_{\mathcal{R}} - [ y ]_{\mathcal{R}} | = | r - y | = r - y \leq x - y = | x - y |$. If $x < l$, $| [ x ]_{\mathcal{R}} - [ y ]_{\mathcal{R}} | = | l - y | = y - l \leq y - x = | x - y |$. Finally, consider the case $x,y \not\in \mathcal{R}$. If $x,y > r$ or $x,y < l$, $| [ x ]_{\mathcal{R}} - [ y ]_{\mathcal{R}} | = 0 \leq | x - y |$. If $x>r, y<l$, $| [ x ]_{\mathcal{R}} - [ y ]_{\mathcal{R}} | = | r - l | = r - l \leq x - y = | x - y |$. The case $x<l, y>r$ follows from a similar argument. \QuodEratDemonstrandum

\subsection{Proof of Lemma \ref{lemma:Ergodic_probability_bound}}
\label{sec:Appendix__Proving_the_probabilistic_error_bound_on_the_steady_state_estimate}

Let $\varianceWrt{ f }{\vect{\mu}} = \frac{1}{2} \sum_{\svA,\svB \in \stateSpace} \bigl( f(\svA) - f(\svB) \bigr)^2 \vectComponent{\mu}{\svA} \vectComponent{\mu}{\svB}$, $(f,g)_{\vect{\mu}} = \sum_{\svA \in \stateSpace} f(\svA) g(\svA) \vectComponent{\mu}{\svA}$ and $\pnorm{ \vect{\mu} }{2,\vect{\nu}} = ( \sum_{\svA \in \stateSpace} \vectComponent{\mu}{\svA}^2 \vectComponent{\nu}{\svA} )^{1/2}$.

\begin{proposition}[\cite{CattiauxGuillin2006}, p.~2]
\label{prop:CattiauxGuillin2006}
On some Polish space $\Omega$, let us consider a conservative (continuous-time) Markov process denoted by $\process{ X(t) }{t \geq 0}$ and with infinitesimal generator $\infinitesimalGenerator{ }$. Let $\vect{\mu}$ be a probability measure on $\Omega$ which is invariant and ergodic with respect to $P_t$. 

Assume that $\vect{\mu}$ satisfies the Poincar\'{e} inequality $\varianceWrt{ f }{\vect{\mu}}  \leq - \kappa (\infinitesimalGenerator{f},f)_{\vect{\mu}}$. Then for all $\theta$ such that $\sup | \theta | = 1$, all $0 < \epsilon \leq 1$ and all $t > 0$, assuming that the initial distribution of $X_s$ is $\nu$,
\begin{align}
&\probabilityBig{ \Bigl| \frac{1}{t} \int_0^t \theta(X(s)) ds - \int \theta d\mu \Bigr| \geq \epsilon } \nonumber \\
&\leq \bigpnorm{ \frac{d\nu}{d\mu} }{2,\mu} \exp{ \Bigl( - \frac{t \epsilon^2}{8 \kappa \varianceWrt{\theta}{\vect{\mu}} } \Bigr) }.
\end{align}
\end{proposition}

\refLemma{lemma:Ergodic_probability_bound} is a direct consequence of \refProposition{prop:CattiauxGuillin2006}. Before we can use \refProposition{prop:CattiauxGuillin2006} to prove \refLemma{lemma:Ergodic_probability_bound}, however, we need to verify all of its assumptions. We will now verify these assumptions for continuous-time, reversible Markov processes with a product form solution. Our method is based on an approach for discrete-time Markov chains \cite{DS91}. 

Define a graph $G=(V,E)$, where $V$ denotes the vertex set in which each vertex corresponds to a state in $\stateSpace$ and $E$ denotes the set of directed edges. An edge $e = ( \svA,\svB )$ is in $E$ if $\phi(e) = \vectComponent{\pi}{\svA} \matrixElement{Q}{\svA}{\svB} = \vectComponent{\pi}{\svB} \matrixElement{Q}{\svB}{\svA} > 0$. Here, $Q$ denotes the generator matrix of $\process{X(t)}{ t \geq 0 } $. For every pair of distinct vertices $\svA, \svB \in \stateSpace$, choose a path $\gamma_{\svA,\svB}$ (along the edges of $G$) from $\svA$ to $\svB$. Paths may have repeated vertices but a given edge appears at most once in a given path. Let $\Gamma$ denote the collection of paths (one for each ordered pair $\svA, \svB$). Irreducibility of $\process{ X(t) }{ t \geq 0 }$ guarantees that such paths exist. For $\gamma_{\svA,\svB} \in \Gamma$ define the path length by $\pnorm{\gamma_{\svA,\svB}}{\phi} = \sum_{ e \in \gamma_{\svA,\svB} } ( 1/\phi(e) )$. Also, let 
\begin{align}
\kappa 
= \max_{e} \sum_{ \{ \gamma_{\svA,\svB} \in \Gamma | e \in \gamma_{\svA,\svB} \} } \pnorm{\gamma_{\svA,\svB}}{\phi} \vectComponent{\pi}{\svA} \vectComponent{\pi}{\svB}
\end{align}
and $f(e) = f(\svB) - f(\svA)$ for $e = (\svA,\svB) \in E$. Then write
\begin{align}
\varianceWrt{ f }{\vect{\pi}}
&= \frac{1}{2} \sum_{\svA,\svB \in \stateSpace} \Bigl( \sum_{ e \in \gamma_{\svA,\svB} } \Bigl( \frac{\phi(e)}{\phi(e)} \Bigr)^{\frac{1}{2}} f(e) \Bigr)^2 \vectComponent{\pi}{\svA} \vectComponent{\pi}{\svB}.
\end{align}
Use the Cauchy-Schwarz inequality $| \transpose{ \vect{x} } \vect{y} |^2 \leq \transpose{\vect{x}} \vect{x} \cdot \transpose{\vect{y}} \vect{y}$ to obtain
\begin{align}
\varianceWrt{ f }{\vect{\pi}}
&\leq \frac{1}{2} \sum_{\svA,\svB \in \stateSpace} \vectComponent{\pi}{\svA} \vectComponent{\pi}{\svB} \Bigl( \sum_{ e \in \gamma_{\svA,\svB} } \frac{1}{\phi(e)} \Bigr) \Bigl( \sum_{ e \in \gamma_{\svA,\svB} } \phi(e) f(e)^2 \Bigr) \nonumber \\
&= \frac{1}{2} \sum_{\svA,\svB \in \stateSpace} \vectComponent{\pi}{\svA} \vectComponent{\pi}{\svB} \pnorm{\gamma_{\svA,\svB}}{\phi} \Bigl( \sum_{ e \in \gamma_{\svA,\svB} } \phi(e) f(e)^2 \Bigr) \nonumber \\
&= \frac{1}{2} \sum_{e \in E} \phi(e) f(e)^2 \sum_{ \{ \gamma_{\svA,\svB} \in \Gamma | e \in \gamma_{\svA,\svB} \} } \pnorm{\gamma_{\svA,\svB}}{\phi} \vectComponent{\pi}{\svA} \vectComponent{\pi}{\svB}. \nonumber
\end{align}
Use the definition of $\kappa$ and the symmetry of $\phi(e)$ to write
\begin{align}
\varianceWrt{f}{\vect{\pi}}
&\leq \frac{\kappa}{2} \sum_{e \in E} \phi(e) f(e)^2 \nonumber \\
&= \frac{\kappa}{2} \sum_{\svA, \svB \in \stateSpace} \vectComponent{\pi}{\svB} \matrixElement{Q}{\svB}{\svA} ( f(\svB)^2 - f(\svB) f(\svA) ) \nonumber \\
&\phantom{=} + \frac{\kappa}{2} \sum_{\svA, \svB \in \stateSpace} \vectComponent{\pi}{\svA} \matrixElement{Q}{\svA}{\svB} ( f(\svA)^2 - f(\svB) f(\svA) ) \nonumber \\
&= \kappa \sum_{\svA, \svB \in \stateSpace} \matrixElement{Q}{\svA}{\svB} ( f(\svA) - f(\svB) ) f(\svA) \vectComponent{\pi}{\svA} \nonumber \\
&= \kappa \sum_{\svA \in \stateSpace} \Bigl( \sum_{\svB \in \stateSpace} \matrixElement{Q}{\svA}{\svB} ( f(\svA) - f(\svB) ) \Bigr) f(\svA) \vectComponent{\pi}{\svA}.
\end{align}
By definition of the infinitesimal generator $\infinitesimalGenerator{ }$, we find that
\begin{align}
(\infinitesimalGenerator{f})(\svA)
&= \lim_{t \rightarrow 0} \frac{1}{t} \Bigl( \sum_{\svB \in \stateSpace} \matrixElement{( \e{t Q} \bigr)}{\svA}{\svB} f(\svB) - f(\svA) \Bigr) \nonumber \\
&= \lim_{t \rightarrow 0} \frac{1}{t} \Bigl( \sum_{\svB \in \stateSpace} \matrixElement{( I + t Q + \bigO{t^2} \bigr)}{\svA}{\svB} f(\svB) - f(\svA) \Bigr) \nonumber \\
&= \sum_{\svB \in \stateSpace} \matrixElement{Q}{\svA}{\svB} f(\svB) 
= \sum_{\svB \in \stateSpace \backslash \{ \svA \} } \matrixElement{Q}{\svA}{\svB} f(\svB) \nonumber + \matrixElement{Q}{\svA}{\svA} f(\svA) \\
&= \sum_{\svB \in \stateSpace \backslash \{ \svA \} } \matrixElement{Q}{\svA}{\svB} f(\svB) \nonumber - \sum_{\svB \in \stateSpace \backslash \{ \svA \} } \matrixElement{Q}{\svA}{\svB} f(\svA) \nonumber \\
&= \sum_{\svB \in \stateSpace } \matrixElement{Q}{\svA}{\svB} ( f(\svB) - f(\svA) ),
\end{align}
after which one can conclude that $\varianceWrt{ f }{\vect{\pi}} \leq -\kappa (\infinitesimalGenerator{f},f)_{\vect{\pi}}$. We also note that when choosing $\theta(X(t)) = \indicator{ X(t) = \svC }$, we have that
\begin{align}
\varianceWrt{ \theta }{\vect{\pi}} 
= \frac{1}{2} \sum_{\svA,\svB \in \stateSpace} \bigl( \indicator{ \svA = \svC } - \indicator{ \svB = \svC } \bigr)^2 \vectComponent{\pi}{\svA} \vectComponent{\pi}{\svB}
\leq \frac{\cardinality{\stateSpace}^2}{2}. \nonumber
\end{align}
Now starting from any state $\svB$, i.e.~the probability distribution with unit mass in state $\svB$, we have for the initial distance
\begin{align}
\bigpnorm{ \frac{d\nu}{d\mu} }{2,\mu}
&= \Bigl( \sum_{\svA \in \stateSpace} \Bigl( \frac{ \vectComponent{\nu}{\svA} }{ \vectComponent{\mu}{\svA} } \Bigr)^2 \vectComponent{\mu}{\svA} \Bigr)^{\frac{1}{2}}
= \frac{1}{ \sqrt{ \vectComponent{\pi}{ \svB } } }
\leq \frac{1}{ \sqrt{ \min_{\svA \in \stateSpace} \vectComponent{\pi}{\svA} } }, \nonumber
\end{align}
since $\vect{\mu} = \vect{\pi}$. Because $\mathcal{R}$ is bounded, $\min_{\svA \in \stateSpace} \vectComponent{\pi}{\svA}$ is bounded from below by some constant $1/\ergodicConstant \in \strictlyPositiveRealNumbers$. \QuodEratDemonstrandum

\subsection{Proof of Lemma \ref{lemma:Mn_is_a_martingale}}
\label{sec:Appendix__Proof_that_Mn_is_a_martingale}

First note that $\iterand{M}{n} \in \iterand{\mathcal{F}}{n}$ and that its expectation is bounded, which can be concluded after writing
\begin{align}
&\expectation{ | \iterand{M}{n} | }
\leq \sum_{j=1}^n \iterand{a}{j} \expectation{ \bigl| \sum_{i=1}^d \iterand{\vectComponent{E}{i}}{j} ( \iterand{\vectComponent{R}{i}}{j-1} - \criticalpoint{\vectComponent{r}{i}} ) \bigr| } \nonumber \\
&\leq \sum_{j=1}^n \iterand{a}{j} \max_{i=1,...,d} \{ \vectComponent{\mathcal{R}^{\max}}{i} - \vectComponent{\mathcal{R}^{\min}}{i} \} \expectation{ \sum_{i=1}^d \bigl| \iterand{\vectComponent{E}{i}}{j} \bigr| }
\end{align}
and then substituting \refEquation{eqn:Error_bias__Upper_bound_in_terms_of_dimension_and_gradient}. Also,
\begin{align}
&\expectation{ \iterand{M}{n} | \iterand{\mathcal{F}}{n-1} }
= \expectation{ \sum_{j=1}^n \iterand{a}{j} \transpose{ \iterand{\vect{E}}{j} } ( \iterand{\vect{R}}{j-1} - \criticalpoint{\vect{r}} ) | \iterand{\mathcal{F}}{n-1} } \nonumber \\
&= \iterand{M}{n-1} + \itstep{n} \expectation{ \transpose{ \iterand{\vect{E}}{n} } ( \iterand{\vect{R}}{n-1} - \criticalpoint{\vect{r}} ) | \iterand{\mathcal{F}}{n-1} }
= \iterand{M}{n-1}, \nonumber
\end{align}
which concludes the proof. \QuodEratDemonstrandum

\end{document}